\documentstyle[epsfig,eqsection]{article}
\newenvironment{pf}{\noindent{\sc Proof}.\enspace}{\rule{2mm}{2mm}\medskip}
\newenvironment{pfn}{\noindent{\sc Proof} \enspace}{\rule{2mm}{2mm}\medskip}
\newtheorem{theorem}{Theorem}
\newtheorem{lemma}{Lemma}
\newtheorem{corollary}{Corollary}

\newtheorem{remark}{Remark}
\newtheorem{definition}{Definition}
\newcommand{\be}{\begin{equation}}
\newcommand{\ee}{\end{equation}}

\newcommand{\th}{\theta}
\newcommand{\Th}{\Theta}
\newcommand{\teta}{\theta}
\newcommand{\Teta}{\Theta}
\newcommand{\om}{\omega}

\newcommand{\lla}{\langle}
\newcommand{\rn}{\rangle}

\newcommand{\ep}{\epsilon}
\newcommand{\al}{\alpha}
\newcommand{\bt}{\beta}

\newcommand{\ov}{\overline}
\newcommand{\wtilde}{\widetilde}
\newcommand{\R}{\mbox{I\hspace{-.15em}R}}
\newcommand{\Rs}{\mbox{ \scriptsize I\hspace{-.15em}R}}
\newcommand{\Z}{\mbox{Z\hspace{-.3em}Z}}

\newcommand{\N}{\mbox{I\hspace{-.15em}N}}

\begin{document}


\title{{\bf Variational construction of  homoclinics 
and chaos in presence of a saddle-saddle equilibrium }}
\author{Massimiliano Berti and Philippe Bolle}
\date{}

\maketitle



We consider autonomous Lagrangian systems with two 
degrees of freedom,
having an hyperbolic equilibrium of saddle-saddle type 
(that is the eingenvalues of the linearized system about the equilibrium are
$ \pm \lambda_1 , \pm \lambda_2 $, $  \lambda_1 , \lambda_2 > 0$). 
We assume that $\lambda_1 > \lambda_2 $ and that 
the system possesses two homoclinic 
orbits. Under a nondegeneracy assumption on the 
homoclinics and under suitable conditions 
on the geometric behaviour of these homoclinics near the equilibrium 
we prove, by variational methods, then they  give rise to
  an  infinite family of multibump homoclinic
solutions and that the topological entropy at the zero energy level is 
positive. A method  to deal also with homoclinics satisfying a weaker 
nondegeneracy condition is developed and it is applied, for
simplicity,  when $\lambda_1 \approx \lambda_2$. An application 
to a  perturbation of a uncoupled system is also given. 
\footnote{msc: 34C37, 58F05, 58E99.}

\section{Introduction}

Let us consider the following Lagrangian system
\be \label{maineq}
-\ddot{q}+\psi (q){\cal J} \dot{q} +Aq = \nabla W (q)
\ee
where $q= (q_1, q_2 ) \in \R^2 $, 
$
{\cal J} = \left(
\begin{array}{cc}
0 & -1 \\
1 & 0
\end{array}
\right)
\ {\rm and} \ A = \left(
\begin{array}{cc}
\lambda_1^2  & 0 \\
0  & \lambda_2^2
\end{array}
\right ).
$
System (\ref{maineq}) can be obtained by the following
Lagrangian   
$$
{\cal L}(q, {\dot q} ) = 
\frac{1}{2} |{\dot q}|^2 + \frac{1}{2}
Aq \cdot q 
- {\dot q} \cdot v(q) -  W(q),
$$
where $v=(v_1,v_2)$ satisfies
\be \label{eq:defv}
 \psi(q) =  \partial_{q_1}{v_2}(q)  -  \partial_{q_2}{v_1}
(q).
\ee
System (\ref{maineq}) admits the energy
$$
{\cal E}( q, \dot{q} ) = \frac{1}{2} |\dot{q}|^2 - \frac{1}{2}
Aq \cdot q + W(q)
$$
as a prime integral. We shall assume
\begin{itemize}
\item (W1)  $ \ W \in C^2 (\R^2 , \R)$, 
$ W(0)=0, \nabla W(0) = 0, D^2 W(0)=0 $; for some 
$0 < \rho_0 < \rho_1 $ ( $\rho_1$ is specified after 
hypothesis $(S2)$)  $D^2 W$ is $L_1$-Lipschitz 
continuous on the ball $B_0 :=
B(0,\rho_0)$ of center $0$ and radius $\rho_0$ 
and $\ov{L}_1$-Lipschitz continuous on $B_1 := B(0,\rho_1)$;
\item ($P1$) 
$\  \psi \in C^1 (\R^2 , \R)$ satisfies $\psi(0)=0$ and is
 $L_2$-Lipschitz continuous
(resp. $\ov{L}_2$-Lipschitz continuous) on $B_0$ (resp. $B_1$); 
 $\nabla \psi$
is $\ov{L}_3$-Lipschitz continuous on $B_1$.
\end{itemize}
By $(P1)$, we can assume
(\ref{eq:defv}), with 
\begin{itemize}
\item $(v1)$  $\ v  \in C^1 (\R^2 , \R^2)$, 
$v(0)=0$, $\nabla v (0)=0$.
\end{itemize}
Under these assumptions $0$ is a 
hyperbolic equilibrium of (\ref{maineq}) and the characteristic 
exponents are two couples of opposite real numbers 
$\pm \lambda_1 , \pm \lambda_2 $. 
In this case 
the equilibrium is called of saddle-saddle type.
We shall assume in the sequel that
$$
\lambda_1 > \lambda_2 >0.  \leqno(S1)
$$
We are interested in a chaotic behaviour of the dynamics
at the zero energy level.
\\[2mm]
\indent 
The only other possibility for a hyperbolic equilibrium
of a Hamiltonian system in a phase space of dimension $4$
is the saddle-focus situation, namely when the 
characteristic exponents
are  $\pm \lambda \pm i  \eta, \ 
\lambda, \eta > 0. $ It would be the case of
system (\ref{maineq}) if 
$ |\psi(0)| \in (|\lambda_1 -\lambda_2| , 
\lambda_1 + \lambda_2)$ 
(note that if $ |\psi(0)| \leq  |\lambda_1 -\lambda_2| $ then
$0$ is still a saddle-saddle equilibrium and that if 
$ |\psi(0)| \geq  \lambda_1 + \lambda_2 $ the equilibrium $0$ is no more
hyperbolic).

The saddle-focus case 
has been investigated by Devaney who showed in \cite{D} 
that, if the system possesses a nondegenerate ( transversal )  homoclinic 
orbit, then it is possible to embed a horseshoe --
and hence a Bernoulli shift -- in the dynamics of the system. This result was
extended by 
Buffoni and S\'er\'e in \cite{BS}, who relaxed
the nondegeneracy condition and proved by 
variational methods the existence of chaos 
at the zero energy level under global assumptions.

These results do not apply in the saddle-saddle case.

The existence of a chaotic dynamics in presence of a saddle-saddle equilibrium
has been studied by Turaev and Shil'nikov \cite{Sh} and more recently by 
Bolotin and Rabinowitz \cite{BR} for a system on a $2$-dimensional torus. 
In this latter paper the existence of homoclinic
orbits is not assumed {\it a priori}, but a simple
geometrical condition is given, which implies that 
the system possesses chaotic trajectories 
either at any small negative energy level
or at  any small positive energy level $\{ {\cal E} = h \}$.
Other  results have been stated in \cite{BS3} for Lagrangian
systems on manifolds.

However, the chaotic trajectories which are obtained 
in \cite{Sh} as well as in \cite{BR} or in \cite{BS3} are not preserved
when the energy vanishes. 
\\[2mm]
\indent
The existence of a 
Bernoulli shift at energy level $ \{ {\cal E } = 0 \} $ 
was studied by Holmes in \cite{H} (see also \cite{Wg}).
He assumed the existence of two nondegenerate homoclinics
and introduced some conditions on the way these homoclinics
approach $0$ which ensure, when
$(S1)$ is satisfied, the existence
of a horseshoe at the zero energy level.
By the structural stability of the horseshoes there results chaos 
also on nearby energy levels $ \{ {\cal E } = h \} $,
see \cite{H}.  
\\[2mm]
\indent
In the present paper we deal as in \cite{H} with
the saddle-saddle case, under assumption $(S1)$.
We give specific 
conditions, called  $(H1-4)$, directly
inspired to the assumptions of Holmes, which imply that 
the system possesses an infinite family of 
multibump homoclinic orbits and of solutions with infinitely many bumps, 
which give rise to a chaotic behaviour at the zero energy level.
Furthermore we improve such results requiring for the homoclinics $\ov{q}$,
$\wtilde{q}$ a nondegeneracy condition weaker than transversality.
Rather than performing this relaxation 
in a general situation, which would 
require quite involved conditions, we restrict ourselves to  the case
when the eigenvalues are close  one to  each other. 
However we underline that 
the method  introduced to deal with degenerate homoclinics 
is could  be 
adapted to a large variety of situations where it 
is difficult or impossible to check the nondegeneracy assumption.

First we shall assume that
\begin{itemize}
\item
$(S2)$ System (\ref{maineq}) has 
2 nondegenerate homoclinics $ {\overline q}, {\widetilde q}$.
``Nondegenerate'' means that the unique solutions of the
linearized equation at (for instance) $\ov{q}$
$$
-\ddot{h} +Ah + \psi(\ov{q}) {\cal J} \dot{h} +
\nabla \psi (\ov{q}) \cdot h  {\cal J} \dot{\ov{q}}
- D^2 W (\ov{q})h=0
$$
that tend to $0$ as $t\to \pm \infty$ are 
$c\dot{\ov{q}}, \ c\in \R$. That means that the
stable and unstable manifolds to $0$ intersect 
transversally at $(\ov{q}(t), \dot{\ov{q}} (t)) $  at the 
zero energy level.
\end{itemize}
\noindent
We can now specify the constant $\rho_1$ in $(W1)$:
$ \rho_1 > \max \{ |\ov{q}|_{\infty} , |\wtilde{q} |_{\infty}  \} 
+ \rho_0 $. 
\\[1mm]
The relaxed nondegeneracy condition is the following 
\begin{itemize}
\item
$(S2')$ System (\ref{maineq}) has 
2 ``topologically nondegenerate'' isolated
homoclinics $ {\overline q}, {\widetilde q} $,
(see definition \ref{eq:topdeg} in subsection \ref{sub:top}).
\end{itemize}
We point out that in some situations such a 
condition can be checked for homoclinics obtained by 
variational methods which are isolated up to time translations, see
for example \cite{A}, \cite{HH}.
\\[1mm]
\indent
In order to get chaotic trajectories 
in the saddle-saddle case it is necessary to postulate the existence 
of (at least) two homoclinic orbits, while only one is necessary for 
the saddle-focus case. 
Even though, there exist systems with  
several  transversal homoclinic orbits 
which do not have a chaotic behaviour.
Consider for example (\ref{maineq}) and assume that:
\be\label{eq:e} 
W(q)= q_1^4 + q_2^4 \qquad {\rm and} \qquad 
\psi(q)=0.
\ee
Then the system reduces to a 
direct product of $1$-dimensional systems.
$(0,0) \in \R^4$ is a saddle-saddle equilibrium with
$4$ transversal homoclinic trajectories but the system
is integrable
(another example of an integrable Hamiltonian system with several 
transversal homoclinic orbits is given in \cite{De}).
Thus additional assumptions are needed for chaotic behaviour.
In order to obtain multibump homoclinics
for system (\ref{maineq}) 
as glued  copies of $\ov{q}$ and $\wtilde{q}$,  
some hypotheses of geometrical  nature on 
$\ov{q}$ and ${\wtilde q}$, similar to the ones given in 
\cite{H}, are required.
\\[5mm]
\indent
The results contained in this paper 
have already been outlined in \cite{BB1}.
In order to describe them  we need some notations.
We shall assume that $\ov{q}(\R) $ and $\wtilde{q} (\R)$ are not included
in $B_0$. 
For $ r \in (0, \rho_0 /2) $ 
  we define $\ov{T} > 0$ by    
$|\ov{q}( \pm \ov{T} ) |= r $ and
$|\ov{q} (t)| < r $ for $ |t| > \ov{T} $. 
We define in the same way $\wtilde{T}$ and we set  
$T= \min \{ \ov{T} , \wtilde{T} \} $.

Call $( {\ov{\al}}_1 , {\ov{\al}}_2 )$ $=$ $(\ov{q}_1(- \ov{T} ), 
\ov{q}_2(- \ov{T}))$, 
$( {\ov{\bt}}_1 , {\ov{\bt}}_2 )$ $=$ $(\ov{q}_1(\ov{T}), \ov{q}_2(\ov{T}))$
the extremal intersection points of 
$\ov{q}( \R ) $ with the circle in ${\R}^2$ of radius $r$;
similarly we introduce 
$( {\wtilde{\al}}_1 ,{\wtilde{\al}}_2 )=
(\wtilde{q}_1(- \wtilde{T} ), \wtilde{q}_2(- \wtilde{T} ))$, 
$( {\wtilde{\bt}}_1 ,{\wtilde{\bt}}_2 )=
(\wtilde{q}_1(\wtilde{T}), \wtilde{q}_2(\wtilde{T}))$. 
Let $\ov{\om}_u$, $\ov{\om}_s$ be defined by
$$ (\ov{\alpha}_1 , \ov{\alpha}_2 ) = ( r \cos 
\ov{\om}_u , r \sin \ov{\om}_u) \quad , \quad  
 (\ov{\beta}_1 , \ov{\beta}_2 ) = ( r \cos 
\ov{\om}_s , r \sin \ov{\om}_s );$$ 
$\wtilde{\om}_u$, $\wtilde{\om}_s$ are defined in the same way.

We set
$\Lambda = (L_1 / \lambda_2^2) + (3 L_2  \lambda_1 / \lambda_2^2 )$,
$\ov{\Lambda}= (\ov{L}_1 / \lambda_2^2) + (3\ov{L}_2 \lambda_1 / \lambda_2^2) +
\max \Big \{|\dot{\ov{q}}|_{\infty},
|\dot{\wtilde{q}}|_{\infty} \Big \} (\ov{L}_3 / \lambda_2^2 )$,
where $L_i$, $\ov{L}_i$ are defined in assumptions $(W1)$, $(P1)$. 
Note that $\Lambda$,  $\ov{\Lambda}$   do not change 
if the equation is modified by a time rescaling $q(t) \to q( \alpha t )$.
\\[1mm]
\indent
In the next conditions $\om_u$ stands for $\ov{\om}_u$ or $ \wtilde{\om}_u$ 
and $\om_s$ for $\ov{\om}_s$ or $\wtilde{\om}_s$.

\begin{itemize}
\item ($H1$) 
$ {\om}_u, {\om}_s
\neq n \pi/2 , \  n \in \Z $, 
$\tan  \om_u   \tan  \om_s  < 0$ 
and ( $ \cos \ov{\om}_u  \cos \wtilde{\om}_u < 0 $ or
$\cos \ov{\om}_s  \cos \wtilde{\om}_s < 0$).

(the above inequalities are satisfied for example if 
$\ov{\om}_u \in ( 0, \pi /2 )$,  
$\ov{\om}_s \in ( 3\pi /2 , 2\pi)$, 
$\wtilde{\om}_u \in ( \pi  , 3 \pi / 2 )$ and
$\wtilde{\om}_s \in ( \pi / 2 , \pi )$);
\item ($H2$) 
$$
\frac{\lambda_2^2}{\lambda_1^2}
\frac{|\alpha_2||\beta_2| +  (15 \lambda_1/4\lambda_2)
\Lambda  r^3}{|\alpha_1||\beta_1|} 
\leq  l(\frac{\lambda_1}{\lambda_2})\min 
\Big ( e^{-2 \frac{\lambda_1 - \lambda_2}{\lambda_2}} ,
\Big ( \frac{C_1 e^{\lambda_2 (T- T_{C_1})}}{18} \Big )^{2
\frac{\lambda_1 - \lambda_2}{\lambda_2}},
\Big ( \frac{C_1^2}{40\ov{\Lambda}r} \Big )^{\frac{\lambda_1 -
\lambda_2}{\lambda_2}} \Big )
$$
where $ l(\nu) = \max_{s \in (0, 1/8)} [(1-s)^2 (1-s/5) / (1+s)^3 ]
s^{\nu  - 1}$ and 
$C_1$ is a constant defined by (\ref{eq:nd}),
which measures the transversality
of the homoclinics: smaller is $C_1$ weaker is the transversality.
$T_{C_1}$ depends only on $C_1$ and $\rho_0$ and it is
defined by (\ref{eq:T_C_1}), section \ref{sec:finite}.
\item
$(H3)$
$$
\frac{\lambda_2^2  }{ \lambda_1^2 }
 \frac{ |\alpha_2||\beta_2|+ (15 \lambda_1 /
4 \lambda_2) \Lambda  r^3}{|\alpha_1||\beta_1|}
\leq l(\frac{\lambda_1}{\lambda_2}) 
\left( \frac{C_1 {\cal M} }{ 36  S_2 +28 \Lambda r^2}
\right)^{(\lambda_1/\lambda_2) -1}
$$
where   ${\cal M} = 
\min \{|\ov{\al}_j|,|\ov{\bt}_j|,|\wtilde{\al}_j |,|\wtilde{\bt}_j| \ ; 
\ j=1,2 \} $ and
$ S_2 = 
\max \{|\ov{\al}_2|,|\ov{\bt}_2|,|\wtilde{\al}_2 |,|\wtilde{\bt}_2| \} $.
\item ($H4$)
$
 \min (| \sin \ov{\om}_{u,s}| , 
| \sin \wtilde{\om}_{u,s} | ) \geq  \sqrt{\frac{\lambda_1}{\lambda_2}
20 \Lambda r} \quad , \quad (12\lambda_1 /\lambda_2 )\Lambda r
\leq C_1.
$ 
\end{itemize}
Roughly speaking the first geometric assumption $(H1)$ means that
the homoclinics $\ov{q}$, $\wtilde{q}$ enter
and leave the origin from different ``quadrants''.
Note that if (\ref{eq:e}) holds system (\ref{maineq}) 
does not satisfy hypothesis (H1).
$(H2-3)$ quantify how small 
$|\tan \om_u \tan \om_s |$ and $r$ must be. Note that 
if the system is linear (that is $W=0$, $\psi =0$) 
in the ball $B(0, \rho_0)$ then
condition $(H4)$ disappears and conditions
$(H2-3)$ are simplified (in $(H2-3)$, $\Lambda =0$).
Moreover if $\lambda_1/ \lambda_2 \to 1 $ then  
$l( \lambda_1 / \lambda_2 ) \to 1$ and the second members in inequalities 
$(H2-3)$ tend to $1$. 
\\[1mm]
\indent
Before stating our first result we introduce some other notations. For
$j=(j_1, \ldots , j_k ) \in \{0 ,1 \}^k$
and for $ \Th = ( \th_1, \ldots, \th_k )$ with 
$\teta_1 <  \ldots < \teta_k $ we define 
$ T_i = \ov{T} $ if $ j_i = 0 $ and $ T_i = \wtilde{T} $ 
if $ j_i = 1 $;
$ d_i =  ( \teta_{i+1}- T_{i+1} )-
(\teta_i + T_i)  $ and $
\ov{d}= \min_{1\leq i \leq k-1} d_i.$

\begin{theorem}\label{thm:main2}
Assume $(W1)$,$(P1)$,$(v1)$,$(S1-2)$ and $(H1-4)$. 
Then there exist $0<D<J$ such that for every $k \in \N $, 
for every sequence $j = (j_1, \ldots , j_k ) \in \{ 0,1 \}^k $ 
there is $\Th=( \teta_1, \ldots ,  \teta_{k} ) \in \R^k $ with 
$ d_i \in ( D, J ) $ for all $ i =1, \ldots, k-1 $ 
and a homoclinic solution of (\ref{maineq}) $x_j$ such that 
\begin{itemize}
\item 
if $j_i =0$ then on the interval $[\teta_i - \ov{T} ,
\teta_i + \ov{T}]$ 
$$
|x_j (t) - \ov{q}(t-\teta_i)| \leq \frac{r}{8}
\min \Big( |\cos \ov{\om}_{u,s}|, |\cos \wtilde{\om}_{u,s}|,
|\sin \ov{\om}_{u,s}|, |\sin \wtilde{\om}_{u,s}| \Big) = \frac{{\cal M}}{8},
$$
\item
if $j_i =1$ then on the interval $[\teta_i - \wtilde{T} ,
\teta_i + \wtilde{T}]$
$$
|x_j (t) - \wtilde{q}(t-\teta_i)| \leq \frac{r}{8}
\min \Big( |\cos \ov{\om}_{u,s}|, |\cos \wtilde{\om}_{u,s}|,
|\sin \ov{\om}_{u,s}|, |\sin \wtilde{\om}_{u,s}| \Big) =
\frac{{\cal M}}{8}, 
$$
\item
Outside $(\cup_{j_i =0} [\theta_i - \ov{T} , 
\theta_i + \ov{T}]) \cup (\cup_{j_i =1} [\theta_i - \wtilde{T},
\theta_i + \wtilde{T}] ) $, $|x_i (t) | \leq 2r$.
\end{itemize}
\end{theorem}

Note that, by theorem \ref{thm:main2}
and  assumption $(H1)$,
two distinct sequences $j=(j_1, \ldots , j_k )$ and 
$j'=(j'_1, \ldots , j'_k )$ give rise to
two distinct homoclinics.

\begin{remark}
{\bf (i)} Since the distance $d_i$ between two consecutive
bumps is bounded by the constant $J$ which is
independent of the number of bumps $k$,  by the Ascoli-Arzel\'a 
theorem  there follows the existence of solutions 
with infinitely many bumps, see {\it theorem  \ref{thm:infin}}.
In particular it implies a lower bound for the topological entropy 
at the zero energy level, 
$h_{top}^0 > \log{2} / (2 \max \{ \ov{T}, \wtilde{T} \} + J)$ and
shows that the system exhibits a chaotic behaviour. 

{\bf (ii)} The fact that $\lambda_1 > \lambda_2 $ is crucial 
to be able to construct multibump homoclinics.

{\bf (iii)} As it will appear in the proof of 
theorem \ref{thm:main2},    smaller are the quantities 
$(\Lambda r / |\cos \om_u \cos \om_s|)  +
|\tan \om_u  \tan \om_s |$,
$|\lambda_1 - \lambda_2 |/\lambda_2$, greater is the 
distance between the bumps.

{\bf (iv)} We do not prove the existence of multibump homoclinics 
in an arbitrary small neighborhood of $\ov{q}, \wtilde{q}$.  
Indeed in \cite{Sh} it is proved that there is a neighborhood
$V$ of $\ov{q}(\R) \cup \wtilde{q}(\R) $ such that  the only homoclinic
solutions contained in $V$ are  $\ov{q}$ and $ \wtilde{q}$.
\end{remark}

Our other results, {\it theorems  
\ref{thm:alm},  \ref{th:pert} and  \ref{thm:rela}}, resp. in 
subsections \ref{sec:almost}, \ref{sub:uncoup} and 
\ref{sub:relth}, 
are variants of theorem \ref{thm:main2} in
special systems or when the homoclinics are degenerate.
\\[2mm]
\indent
The multibump homoclinic solutions of (\ref{maineq})
will be obtained as critical points
of  the following action functional, which is well defined
by $(W1)$ and $(v1)$ on $E=W^{1,2}(\R, \R^2)$: 
\be \label{def:f}
f( q) =  \int_{\Rs} \frac{1}{2} |\dot{q}|^2 +
\frac{1}{2} Aq \cdot q  -\dot q \cdot v(q) -  W(q).
\ee
The idea of the proofs goes as follows.

A `` pseudo-critical''  manifold for $f$, 
$Z_k$ = $\{ Q_{\Th} \ | \ \Th \in {\R}^k, \ \th_1 < \ldots < \th_k  \} $ 
is constructed by  gluing together  translates  of the homoclinics 
$\ov{q}( \cdot - \th_i ) $ and 
$\wtilde{q}( \cdot - \th_j )$,
see section (\ref{sec:bound}) and (\ref{sec:nat}).  
Then we show  that, 
when the bumps are sufficiently separated, that is when
\be \label{eq:sep}
 \min_i  (\th_{i+1}-\th_i ) > \ov{D}
\ee
a shadowing type lemma
enables to construct  
immersions 
${\cal I}_k : \ 
M_k = \{ \Th \in {\R}^k \ | \ \min_i  (\th_{i+1}-\th_i ) > \ov{D} \} 
\to E $ with $ {\cal I}_k ( M_k ) \approx Z_k $   
such that the critical points of  $g(\Th)= 
f ({\cal I}_k ( \Th )) $     
gives rise to a $k$-bump homoclinic solutions.
The geometric properties $(H1-4)$
of the homoclinics $\ov{q}$ and $\wtilde{q}$
ensure the existence of critical points of $g(\Th)$
satisfying (\ref{eq:sep}). We point out that
 $g(\Th)$ does not possess critical points  when
$\min_i ( \th_{i+1} - \th_i ) \to + \infty $; 
therefore  we need to estimate 
carefully the minimal distance $\ov{D}$ for which we 
obtain the immersions ${\cal I}_k$. 
This is done in section \ref{sec:finite}. 
\\[1mm]
\indent
For the sake of clarity we perform all the detailed computations 
for a system with $2$ degrees of freedom, but 
the same method can be adapted  also to study systems in dimension $n$,
(see remark \ref{rem:dimn}) where the analytical technics  based on 
the study of Poincar\'e sections are more difficult.
\\[1mm]
\indent
The paper is organized as follows. 
In section \ref{sec:finite} 
we perform the finite dimensional reduction for the functional $f$ and 
we prove {\it thm. \ref{thm:main2}}.
In section \ref{sec:exam} we give 
examples of applications of theorem  
\ref{thm:main2} when the eigenvalues are near one each other 
({\it thm. \ref{thm:alm}})
and for 
a system which is 
a perturbation of $2$-uncoupled Duffing equations 
({\it thm.  \ref{th:pert}}). 
In section \ref{sec:deg} it is shown that
in the case $\lambda_1 \approx \lambda_2$ 
the transversality condition can be weakened 
assuming the topological nondegeneracy $(S2')$ 
({\it thm. \ref{thm:rela}}).   
Finally in section \ref{sec:dyn}
we show why the above theorems imply a chaotic dynamics
({\it thm. \ref{thm:infin}}).

\section{Finite dimensional reduction}\label{sec:finite}

We shall use the following Banach spaces:
\begin{itemize}
\item 
$Y = W^{1,\infty} (\R , \R^2)$ endowed with norm
$
||y|| = \max \Big( |y|_{\infty} ,\frac{1}{\lambda_2} |\dot{y}|_{\infty}
\Big)
$
where
$|y|_{\infty} = \sup_{t\in \Rs} |y(t)|$.
\item
$E= W^{1,2}(\R , \R^2)$ endowed with scalar product
$ ( x, y ) =  \sum_{j=1}^2  \int_{\Rs} 
{\dot x}_j {\dot y}_j  + \lambda_j^2 x_j y_j $ 
and associated norm $| \cdot |_E $.
\item
$X = \{ h \in Y  \ | \  e^{\lambda_2 |t| } |h(t)|, 
\ e^{\lambda_2 |t|} |\dot{h}(t)|  \in L^{\infty} \}$.
\end{itemize}
Since the equilibrium $0$ is hyperbolic of smaller positive
 characteristic exponent
$\lambda_2$  by standard results
(see also lemma \ref{lem:homo}) any homoclinic solution to $0$ 
of (\ref{maineq})  belongs to $X$.

We have  $X \subset Y \cap E$. 
For $A \subset X$ we shall use the notation
$$
A^{\bot} = \{ y \in Y \ | \
(a,y) =0, \quad \forall a \in A  \}. 
$$
Note that, by the exponential decay of the elements of $X$,
$A^{\bot}$ is well defined and it is a closed subspace of $Y$.
\\[1mm]
\indent
We define the operator $S : Y \to Y$ by 
$$
S(y) = y - L_A (\nabla W (y) - \psi(y) {\cal J} {\dot y})  
$$
where $L_A$ is the linear operator  which assigns to $h$ 
the unique solution $z = L_{A} h $ of 
$$ 
- {\ddot z} + A z = h \qquad {\rm with} \qquad
\lim_{|t| \to \infty } z (t) = 0.
$$
An explicit definition of  $L_A$ is 
\be \label{eq:L_A}
(L_{A} h)(t) = 
\frac{(\sqrt{A})^{-1}}{2} 
\int_{-\infty}^{+\infty} e^{ - |t-s| \sqrt{A} } h(s) ds \  ,  \ 
(\frac{d}{dt} L_A h)(t) = \frac{1}{2} 
\int_{-\infty}^{+\infty} sgn(s-t) \,  e^{ - |t-s| \sqrt{A} } h(s) ds     
\ee
By (\ref{eq:L_A})
it is easy to see that, for all $x \in Y$,
\be \label{eq:estL_A}
||L_A x || \leq \frac{1}{\lambda_2^2} |x|_{\infty}
\leq \frac{1}{\lambda_2^2} ||x||.
\ee
By (\ref{eq:estL_A}) and $(W1)-(P1)$  
we see that the operator $S$ is $C^1$ on $Y$. 
We can also get the straightforward estimate
\be \label{eq:estnodS}
\forall \ ||y|| < \rho_1 \qquad  
||dS(y)h|| \leq \Big( 1+\ov{\Lambda} ||y|| \Big) ||h||.
\ee 
Note also that $S(E \cap Y) \subset E \cap Y$ and that
 for all $q, x \in E \cap Y$
$$ ( S(q), x ) = df(q)[x].$$
If  $S(q)=0$ and $q \in E \cap Y$ then $q$ is a homoclinic 
solution  to system (\ref{maineq}). We can say a 
little bit  better. 
\begin{lemma} \label{lem:one}
Assume that $q \in Y$ satisfies $S(q) =0$ and that 
$\limsup_{|t| \to +\infty} \max(|q(t) |,
|\dot{q}(t)| / \lambda_2 )  < \min (2/\Lambda , \rho_0) $. 
Then $q$ is a homoclinic solution to
(\ref{maineq}).
\end{lemma}
\begin{pf}
Let $m(t) = \max(|q(t) |,
|\dot{q}(t)| / \lambda_2 )$ and
$c= \limsup_{|t| \to \infty } m(t)$. We assume that 
$c < \min (2/\Lambda , \rho_0)$ and we want to prove that 
$c=0$.  Provided $ m(t) \leq \rho_0 $ we have
$$
\Big| \nabla W (q(t)) - \psi (q(t)) {\cal J} \dot{q}(t) \Big|  \leq 
\frac{1}{2}\lambda_2^2  \Lambda (m(t))^2.
$$
Now easy estimates in the expression of $L_A$ show 
that, if $h \in Y$, then
$$
\limsup_{|t| \to \infty} \max \Big( \Big| L_A h (t) \Big| , \frac{1}{\lambda_2}
\Big| \frac{d}{dt}{L_A h } (t) \Big| \Big) \leq \frac{1}{\lambda_2^2} 
\limsup_{|t| \to \infty} |h(t)|. 
$$
Therefore, since $S(q) =0$, we get
$ c \leq  \Lambda c^2 / 2 $, which implies
$c=0$ by our assumption.
\end{pf}

\begin{remark} \label{rem:S'(q)}
If $q$ is a homoclinic solution to (\ref{maineq}) then, by the
characteristic exponents  of the equilibrium, all $y\in Y$ which 
satisfies $dS (q)\cdot  y =0$  belongs to $X$.
So the nondegeneracy condition $(S2)$
amounts to assuming  that Ker $ dS(\ov{q})$ is spanned 
by $\dot{\ov{q}}$, where $dS(\ov{q})$ is regarded as a linear
operator from $Y$ to $Y$.
Moreover $dS(q)$ has the form $Id + K$, where $K$ is a 
compact operator on $Y$. 
 In addition 
$dS(\ov{q}) (Y) \subset {\ov Y}'$, where ${\ov Y}' = \dot{\ov{q}}^{\bot}$. 
Hence  $dS(\ov{q})$ is a 
linear automorphism of ${\ov Y}'$.
\end{remark}

We now introduce another supplemetary space $\ov{Y}''$ to $\dot{\ov{q}}$. 
The introduction of the above norm $|| \cdot ||$ and  $\ov{Y}''$ 
instead the more natural $H^1$-norm and $\dot{\ov{q}}^{\bot}$ 
is motivated by the fact that this choise allows to obtain 
better estimates in hypotheses $(H2-4)$.

Consider $\ov{t}$ such that $|\dot{\ov{q}}(t)|$ attains its maximum 
at $\ov{t}$.
Let $\ov{\tau}$ be some positive real number such that 
$|\dot{\ov{q}} (t)| \geq 3  |\dot{\ov{q}} (\ov{t})| / 4  $ on the interval 
$\ov{J} = (\ov{t} - \ov{\tau} , \ov{t} + \ov{\tau})$. 
Let 
$$
 \ov{a}_0 = L_A ( \dot{\ov{q}} {\chi}_{\ov{J}}), 
$$
where $\chi_{\ov{J}} $ is the characteristic function of the
interval $\ov{J}$.
By the expression of $L_A$ (\ref{eq:L_A}), we see that 
${\ov a}_0   \in X $.
We define
$$\ov{Y}'' = {\ov{a}}_0^{\bot} = 
\Big\{ h \in Y \  \Big| \  \int_{\ov{t}- \ov{\tau}}^{\ov{t}+ \ov{\tau} }
 \dot{\ov{q}}(s)
\cdot h(s) \; ds =0 \Big\} .$$ 

$\ov{Y}''$ is a supplementary to $ \dot{\ov q} $ and hence
by remark \ref{rem:S'(q)} there exist 
a positive constant $C_0$ such that for all $h \in \ov{Y}''$
$$
\min_{\mu \in \Rs}|| dS({\ov q}) h - \mu \ov{a}_0 || \geq C_0 ||h|| 
$$
(note that, due to the fact that $dS(\ov{q}) = Id +\, compact$, $C_0 \leq 1$).
This  implies that  
\be \label{eq:nd}
\max \Big( ||dS({\ov q}) h - \mu \ov{a}_0 ||, 
R |(h, \ov{a}_0 )| \Big) \geq C_1 ||h||, \quad 
\forall (h, \mu ) \in Y \times \R, 
\ee
where the constant $C_1 \leq 1$ and $ C_1 \to C_0 $ 
as $ R \to +\infty $. In the sequel 
we will fix $R$ and assume that (\ref{eq:nd}) holds
(we can choose $C_1$ as close to $C_0$ as desired).

Now let $ \ov{a} = \alpha \ov{a}_0 $, where $ \alpha > 0 $ is chosen such 
that $ ||\ov{a}|| = C_1 + 1 + \ov{\Lambda} ||\ov{q}|| $.

It is easy to see by (\ref{eq:estnodS}) and (\ref{eq:nd}) that 
\be \label{eq:tra2}
\max \Big( ||dS(\ov{q}) h - \mu \ov{a}||, R|(h, \ov{a} )| \Big)  \geq 
C_1 \max ( ||h||, |\mu| ), \qquad \forall (h, \mu) \in Y \times \R  .
\ee 
We shall assume that also for  $ \wtilde{q} $ are defined the   
corresponding quantities $ \wtilde{t} $, 
$\wtilde{\tau}$, $\wtilde{a}$ and that 
condition (\ref{eq:tra2}) holds.
In the sequel we will also assume that 
$ \max \{ \ov{\tau}, \wtilde{\tau} \}  < T $. 

Now we define $T_{C_1}$. Let $\ov{T}_{C_1}$ the smallest positive 
time such that
\be \label{eq:T_C_1}
\forall t \in \R \backslash [-\ov{T}_{C_1}, \ov{T}_{C_1}] \quad
|\ov{q}(t)| \leq \rho_0  
\quad {\rm and} \quad 8 \Lambda \max \Big( |\ov{q}(t)|
 , \frac{|\dot{\ov{q}}(t)|}{\lambda_2}\Big) \leq C_1 .
\ee
We can define in the same way $\wtilde{T}_{C_1}$, and we set
$T_{C_1} = \max ( \ov{T}_{C_1} , \wtilde{T}_{C_1} )$.

The reason for this definition will appear in the 
proof of lemma \ref{lem:deriv}. 
It is easy to see that, if  ($H4$) is satisfied, then 
by lemma \ref{lem:homo} we have $\ov{T}_{C_1} \leq \ov{T} $,
$\wtilde{T}_{C_1} \leq \wtilde{T}$.

\subsection{Boundary value problems}\label{sec:bound}

The aim of this section is to show how solutions of 
the non-linear system (\ref{maineq}) 
are approximated by solution of the 
linear one $ - {\ddot q} + A q  =  0 $  in a
sufficiently small neighborhood of the origin 
$B_r = \{ q \in {\R}^2 \ | \ |q| \leq r \}$.

First we consider the linear case. 
The solution 
$ q_{d,L} (t) : [0,d] \to \R^2 $ of the linear system
$ - {\ddot q} + A q  =  0 $ with boundary conditions
$ q_{d,L} ( 0 ) = \bt $ and $q_{d,L} ( d ) = \al $ is given by
\be{\label{eq:linear}} 
(q_{d,L})_j (t)  = 
\frac{ \bt_j \sinh ( \lambda_j (d-t) ) + 
\al_j \sinh ( \lambda_j  t ) }{ \sinh ( \lambda_j d ) }, 
\quad j=1,2, 
\ee
whereas the  solutions
$
q^{+}_{h,L} : [ 0, +\infty) \to B_r $  (resp.      
$q^{-}_{h,L} : (-\infty, d ] \to B_r   
$)    
of the linear system $ - {\ddot q} + A q = 0 $
such that $ \lim_{t \to + \infty} q^{+}_{h,L} (t) = 0 $ (resp.
$ \lim_{t \to - \infty} q^{-}_{h,L} (t) = 0 $) and 
$q^{+}_{h,L} (0) = \bt $ (resp. $q^{-}_{h,L} (d) = \alpha$) 
are given by
\be\label{eq:homlin}
q^{+}_{h,L}(t)=  e^{- t \sqrt{A}} \bt,  \qquad     
q^{-}_{h,L}(t) = e^{ (t - d) \sqrt{A} } \alpha.   
\ee
We define
\be\label{eq:deltal}
\Delta_L^d (\beta , \alpha )= \max \Big\{ 
\Big| {\dot q}_{d,L}(0) - {\dot q}^{+}_{h,L}(0) \Big|, 
\Big|{\dot q}_{d,L}(d) - {\dot q}^{-}_{h,L}(d) \Big| \Big\}. 
\ee
By (\ref{eq:linear}) and (\ref{eq:homlin}) we can compute
$$
\displaystyle
({\dot q}_{d,L}(0) - {\dot q}^{+}_{h,L}(0))_j =  
 \frac{ \lambda_j (\al_j - \beta_j e^{-\lambda_j d})}
{ \sinh ( \lambda_j d ) } \quad , \quad   
({\dot q}_{d,L}(d) - {\dot q}^{-}_{h,L}(d))_j =  \frac{ \lambda_j 
(- \bt_j + \alpha_j e^{-\lambda_j d})}
{ \sinh ( \lambda_j d  ) }. 
$$
We shall always assume that $d\geq 2/\lambda_2$. Setting
$S_j = \max (|\alpha_j|, |\beta_j|) $, we deduce that
\be\label{eq:delta}
\Delta_L^d \leq 
2 \sqrt{2} \Big( \frac{1+e^{-2}}{1-e^{-4}} \Big)
\max_{j=1,2} \Big\{  \lambda_j S_j \exp -(\lambda_j d ) \Big\}. 
\ee

We now consider the analogous solutions of the non-linear system.
Since $0$ is a hyperbolic equilibrium the existence of the local  
stable and unstable manifolds is standard. 
The following lemma would follow from that but 
we prove it directly by a fixed point argument 
because we need some  explicit estimates. 

\begin{lemma} \label{lem:homo}
For all $ 0 < r < r_0 $ with $r_0 = \min 
(1/6 \Lambda , \rho_0 /2 ) $, for all
$\alpha, \bt \in {\R}^2 $ with $|\alpha|=| \bt | = r $
there exist unique trajectories of (\ref{maineq})  
$$
q^{+}_h : [ 0, +\infty) \to B_r \quad {\rm and } \quad     
q^{-}_h : (-\infty, d ] \to B_r   
$$    
such that $ \lim_{t \to + \infty} q^{+}_h (t) = 0 $,
$ \lim_{t \to - \infty} q^{-}_h (t) = 0 $ and 
$q^{+}_h (0) = \bt $, $q^{-}_h (d) = \alpha$.   

Moreover for all $t $ we have that
$$
| q^+_h (t) - q^+_{h,L} (t)| \leq  \frac{2}{7} r^2 \Lambda e^{-\lambda_2 t }
\quad , \quad 
|{\dot q}^+_h (t) - {\dot q}^+_{h,L} (t)| \leq  \frac{2}{7}
\lambda_2  r^2 \Lambda e^{-\lambda_2 t},
$$
$$
|q^{+}_h ( t)| \leq  \frac{22}{21}r e^{ - \lambda_2 t } \quad , \quad
|\dot{q}^{+}_h ( t)| \leq  \frac{22}{21} \lambda_2 r e^{ - \lambda_2 t }
$$
and the corresponding estimates for $q^-_h$.
\end{lemma} 
The proof of lemma \ref{lem:homo} is given in the appendix.

Now we give a lemma, 
which will be used to glue 
together consecutive bumps, on the existence and uniqueness of orbits  
connecting two points $\al , \bt $
in a neighborhood $B_r$ of $0$. 
Such kind of  lemma is  certainly  not new being deeply
related with the $\lambda$-lemma. However,
since we want to obtain specific estimates, 
we will give a proof based on a fixed point argument. 

\begin{lemma}\label{lem:link1}

For all $ 0 < r < r_1 $ with $ r_1 = 
\min ( 1 / 10 \Lambda , \rho_0 /2) $, for all
$\al, \bt \in {\R}^2 $ with $ | \bt | = | \al | = r $, for all 
$ d > 2 / \lambda_2 $
there exists a unique trajectory of 
(\ref{maineq}) $ q_d (t) $ 
such that $q_d (0) = \bt $, $ q_d ( d ) = \al $
and $q_d ([0,d]) \subset B(0,2r)$.

Moreover the following estimate holds:
\be\label{eq:deltaQ}
\left|\Big({\dot q}_d (0) - {\dot q}^{+}_h (0) \Big) - 
 \Big( {\dot q}_{d,L} (0) - {\dot q}^{+}_{h,L} (0) \Big) \right|,
\left| \Big( {\dot q}_d (d) - {\dot q}^{-}_h (d) \Big) -
\Big( {\dot q}_{d,L} (d) - {\dot q}^{-}_{h,L} (d) \Big)  \right|
 < 5  \lambda_2 r^2 \Lambda  e^{-\lambda_2 d} 
\ee 
\be \label{eq:deltaQ2}
|q_d(t) - q^+_h (t)| < \frac{7}{5} r e^{\lambda_2 (t-d)} \quad , \quad 
|q_d(t) - q^-_h (t)| < \frac{7}{5} r e^{-\lambda_2 t}.
\ee
\be \label{eq:delta-Q2}
|\dot{q}_d(t) - \dot{q}^+_h (t)| < \frac{14}{5} \lambda_1 r 
e^{\lambda_2 (t-d)} \quad , \quad 
|\dot{q}_d(t) - \dot{q}^-_h (t)| < \frac{14}{5}\lambda_1 r 
e^{-\lambda_2 t}.
\ee
\end{lemma}

The proof of lemma \ref{lem:link1} is given in the appendix.
In the sequel we will call also 
$\gamma( \bt, \al, d ) = q_d $ the connecting
solution given by lemma \ref{lem:link1}  and 
$$ 
e(\bt, \al, d)=
\frac{1}{2} \int_0^d {\dot q}^2_d (t) + A q_d (t) \cdot q_d (t) dt 
- \int_0^d {\dot q}_d \cdot v(q_d) dt - \int_0^d W(q_d(t)) dt,
$$ 
the value of the action on the solution $q_d = \gamma ( \bt , \al , d )$.

\subsection{Natural constraint}\label{sec:nat}

As stated in the introduction our existence results 
are obtained by means of a finite dimensional 
reduction according to the following
definition:

\begin{definition}
Let $M$ be a manifold. An immersion  ${\cal I} : M \to  Y$ 
such that ${\cal I}(M) \subset E$ 
is said to induce a natural constraint for the functional $f$
if  
$$
\forall x \in M,  \    \ d(f \circ {\cal I}) (x) = 0 \ 
implies \ that  \ df ({\cal I}(x)) = 0.
$$
\end{definition}

In the sequel for $ j =(j_1, \ldots , j_k ) \in \{0 ,1 \}^k $
and for $\Teta = (\teta_1 , \ldots , \teta_k ) \in \R^k$, with 
$ \th_1 < \ldots < \th_k $ we will use the following notations: 
$$ J_i = \th_i + ( \ov{t} - \ov{\tau}, \ov{t} + \ov{\tau} ) 
\quad {\rm if} \quad 
j_i = 0 \qquad {\rm  and}  \quad
J_i = \th_i + (\wtilde{t} - \wtilde{\tau}, 
\wtilde{t} + \wtilde{\tau} ) \quad {\rm  if} \quad 
j_i = 1.
$$
$$ 
u_i = \th_i - T_i \quad {\rm and } \quad s_i = \th_i + T_i,
$$ 
$$
d_i =  ( \teta_{i+1}- T_{i+1} )-
(\teta_i + T_i) = u_{i+1} - s_i 
\qquad 
{\rm and}
\qquad
\ov{d}= \min_{1\leq i \leq k-1} d_i.
$$
For simplicity the dependence on $\Th=
(\th_1, \ldots , \th_k) $ of $s_i$, $u_i$, $d_i $ and $\ov{d}$ will
remain implicit. 
\\[1mm]
\indent
Fix $j = (j_1, \ldots , j_k ) \in \{0 ,1 \}^k $. 
Our aim is to
prove the existence of a $k$-bump homoclinic 
associated to $j$.
We now define the ``pseudo-critical manifold''.

Consider the $k$ parameter family of continuous 
functions $Q_{\Th}$ defined in the following way:
$$
\begin{array}{rcl}
Q_{\Th}= 
\end{array}
\left\{
\begin{array}{rcl}
& & 
Q^1 ( t ) \; {\rm if} \; t \in  ( - \infty , s_1 ], \\
& & 
\gamma ( Q^1(s_1) , Q^2 (u_2) , d_1 )( \cdot - s_1 ) 
  \; {\rm if} \; t \in [ s_1  , u_2 ],  \\
& & \ldots \\
& & Q^i ( t ) \; {\rm if} \; t \in  [ u_i , s_i ], \\
& & 
\gamma ( Q^i ( s_i )  , Q^{i+1} ( u_{i+1} ) , d_i )( \cdot - s_i ) 
\; {\rm if} \; t \in [ s_i  , u_{i+1} ],  \\
& & \ldots \\
& &
Q^k( t) \; {\rm if} \; t \in [ u_k , + \infty )
\end{array}\right.
$$
where
$$
\begin{array}{rcl}
Q^i ( t ) = 
\end{array}
\left \{ \begin{array}{rcl}
& & 
\ov{q} ( \cdot - \th_i  ) \quad {\rm if} \quad j_i = 0 \\
& &
\wtilde{q} (\cdot - \th_i) \quad {\rm if} \quad j_i = 1 .
\end{array}\right.
$$
We recall that $\gamma$ is defined in lemma \ref{lem:link1}. 
The $k$-dimensional manifold:
$$Z_k = \{ Q_{\Th},
 \ \Th \in \R^k, \  \ov{d} >  2 / \lambda_2  \} $$
is a  $k$-dimensional ``pseudo-critical'' manifold for $f$.
This means that
$|| S ( Q_{\Th} )|| \to 0$ as $ \ov{d} \to + \infty$.
We will give in lemma \ref{lem:grad} a more precise estimate. 
\\[1mm]
\indent
We now show, following \cite{BB}, how to build the 
immersions ${\cal I}_k$.
For $(h, \mu) \in Y \times \R^k $ we denote
$||(h, \mu)|| = \max (||h|| , |\mu_1|, \ldots , |\mu_k | ).$
Let us define the function
$$
H \ : \  \R^k \times Y \times {\R}^k \to Y \times \R^k   
$$
with components $H_1 \in Y $ and $H_2 \in \R^k $ given by:
\begin{eqnarray*}
H_1 (\th_1 , \ldots , \th_k , h, \mu_1 , \ldots , \mu_k ) &=&
S ( h) - 
\sum_{i=1}^k   \mu_i a_i  \\    
H_2 ( \th_1, \ldots , \th_k ,h, \mu_1, \ldots , \mu_k ) &=&
\Big( R(h -Q_{\Teta}, a_1  ), \ldots , 
R(h-Q_{\Teta},  a_k ) \Big) 
\end{eqnarray*}
where
$$
\begin{array}{rcl}
a_i = 
\end{array}
\left \{
\begin{array}{l}
\
\ov{a} ( \cdot - \teta_i )  \quad  {\rm if} \quad j_i=0 \\
\
\wtilde{a} ( \cdot - \th_i ) \quad {\rm if} \quad j_i=1.
\end{array}\right.
$$
Note that $H$ is a $C^1$ function of $(\Teta , h, \mu)$
($Q_{\Teta}$ is not a $C^1$ function of $\Th$
but the function  
$\Teta \mapsto (Q_{\Teta}, a_i) = 
\alpha \int_{J_i} \dot{Q}^i (t)\cdot Q_{\Teta}(t) dt
=\alpha \int_{J_i} \dot{Q}^i (t)\cdot Q^i(t) dt $ is constant.

Consider the partial derivative of $H$,
$
 \partial H / \partial (h,  \mu) , 
$
evaluated at
$( \Th , Q_{\Teta} , \mu )$.
It is the linear operator of $ Y \times \R^k $ 
given by:
\begin{eqnarray*}
\frac{\partial H_1}{\partial 
(h , \mu )}_{|( \Th , Q_{\Th} , \mu )}
[ x , \eta_1 , \ldots , \eta_k ] &=&
dS( Q_{\Th} )x - 
\sum_{i=1}^k \eta_i a_i  \\      
\frac{\partial H_2}{\partial (h , \mu)}_{|(  \Th , Q_{\Th} , \mu )}
[ x , \eta_1 , \ldots ,  \eta_k ] &=&
\Big( R ( x, a_1 ), \ldots ,R ( x, a_k ) \Big) .
\end{eqnarray*}
Note that it is of the form $Id$ $+$ $Compact$ and that it is independent 
of $\mu$ ( and so  we shall omit to write $\mu$ ). 

Following \cite{BB} there results that, 
provided $\ov{d}$ is great enough,
$ \displaystyle \frac{ \partial H}{ \partial ( h, \mu )}_{|(\Th , Q_{\Th})}$
is invertible also on the pseudo-critical manifold $Z_k$ 
and the norm of the inverse satisfies a uniform bound.
As said in the introduction we need in this case a 
specific estimate on $\ov{d}$.

\begin{lemma}\label{lem:deriv}
  
Let $D_1 =\max ( 2/\lambda_2 , 2 (\ln (18/C_1)/ \lambda_2) - 2(T-T_{C_1}))$,
 and assume that $\ov{d} \geq D_1$ and that (H4) holds. 
Then, for all $x \in Y$, for all $\eta =(\eta_1, \ldots , \eta_k ) \in \R^k $,
we have
$$
\Big|\Big|
\frac{\partial H}{\partial (h , \mu )}_{|( \Th, Q_{\Th} )}[x, \eta ] 
\Big|\Big| \geq \frac{C_1}{2} ||(x, \eta)||
$$
{\it i.e.} :
$$
\Big|\Big|  \Big( dS(Q_{\Teta}) x- \sum_{j=1}^k  \eta_j a_j, R (x, a_1) ,
\ldots , R(x,a_k) \Big) \Big|\Big| \geq \frac{C_1}{2} ||(x, \eta_1, \ldots
, \eta_k)||.
$$
\end{lemma}
 
\begin{pf}
We set $m_i = [(\teta_{i+1} -T_{i+1} )+(\teta_i +T_i)] /2$, 
$I_1 = (-\infty , m_1]$, $I_i =
[m_{i-1}, m_i ] $ for $2\leq i \leq k-1$, 
$I_k = [ m_{k-1}, +\infty)$ and  
$||x||_i = \max ( \sup_{t\in I_i} |x(t)|,
\sup_{t \in I_i} |\dot{x} (t) | /\lambda_2 )$.

We define also the compact operators $K, \ov{K}_i , K_i$:
$Y \to Y$ by 
$$
(K x)(t)=L_A {\cal R}x, \qquad (\ov{K}_i  x)(t)=L_A \ov{\cal R}_i x,
\qquad (K_i x)(t)=L_A {\cal R}_i x,
$$
where
$$
{\cal R} x (s) = D^2 W (Q_{\Teta}  (s)) x(s) -
 \psi ( Q_{\Teta} (s)) {\cal J} {\dot x}(s) -
\nabla{\psi} ( Q_{\Teta} (s)) \cdot x {\cal  J} {\dot Q}_{\Teta}(s),
$$ 
$$
\ov{{\cal R}}^ix (s)= D^2 W (Q^i  (s)) x(s) - 
\psi ( Q^i (s)) {\cal J} {\dot x}(s) -
\nabla{\psi} (Q^i(s)) \cdot x {\cal  J}  {\dot Q}^i (s),
$$
${\cal R}_i x (s) = {\cal R}x (s)$ for 
$ s \in {\cal J}_i := [ \teta_i - T^i_{C_1} , \teta_i + T^i_{C_1}]$,
and ${\cal R}_i x (s) = 0 $ on $\R  \backslash {\cal J}_i   $. Here 
we use the notation $T_{C_1}^i = \ov{T}_{C_1}$ if $j_i =0$,
and $T^i_{C_1} = \wtilde{T}_{C_1}$ if $j_i =1$.

For all $x \in Y$ we can write
$
dS(Q_{\Teta} )x = x - K x.
$
We shall prove that 
\be \label{eq:estK_i}
\forall i \, ||\ov{K}_i x - K_i x || \leq (C_1/8)||x|| \ , \quad
||K x - \sum_{i=1}^k K_i x || \leq (C_1/8)||x||.
\ee
We first derive the lemma from (\ref{eq:estK_i}).
Let $(x, \eta ) \in Y \times {\R}^k $; set 
$$ {\cal S} = 
 \Big| \Big|  \Big( dS(Q_{\Teta}) x- \sum_{j=1}^k  \eta_j a_j, R (x, a_1) ,
\ldots , R(x,a_k) \Big) \Big|\Big|  
=  \max \Big\{ || x - Kx - \sum_{j=1}^k 
 \eta_j a_j ||_i, R |( a_i , x ) | \, ; \, 1 \leq i \leq k   \Big\}.
$$
By the second inequality of (\ref{eq:estK_i})
\be \label{eq:cat}
{\cal S} \geq   {\max}_i \Big\{ \max  \Big\{ || x - K_i x - \eta_i a_i ||_i, 
R |( a_i , x ) | \Big\}  - 
\sum_{j \neq i } \Big( || K_j x  +  \eta_j  a_j||_i \Big)
\Big \} - \frac{C_1}{8} ||x||.
\ee

We now define $ z_i \in Y $ by $ z_i = x $ on
$ I_i $ and by
$ -\ddot{z}_i + Az_i = 0 $  on $ (-\infty, m_{i-1}) 
\cup (m_i , +\infty ) $  $  \Big( (m_i , +\infty) $  if 
$  i=1, $ $  (-\infty, m_{k-1}) $  if $  i=k \Big),$ 
$ \lim_{t \to \pm \infty} z_i (t)=0.$
By the definition of $ z_i $ and $ a_i $ it is easy to see that
$(z_i , a_i)= (x, a_i)$,
$||z_i|| = ||x||_i$ and that 
$||z_i -  K_i z_i - \eta_i a_i|| = ||x -  K_i x 
- \eta_i a_i||_i $.
Moreover setting $ M_i = 
\max \Big( ||z_i - K_i z_i -  \eta_i a_i ||, R |(z_i , a_i)| \Big)$ 
and $M= \max_i M_i$, we know by (\ref{eq:tra2}) and 
(\ref{eq:estK_i})  that for all $i$ 
\be \label{eq:Mi}
M_i \geq   C_1 \max (||z_i||, |\eta_i| ) - ||K_i z_i - \ov{K}_i z_i||
\geq (7C_1 / 8) \max (||z_i||, |\eta_i| ). 
\ee
As a consequence
$ M \geq (7C_1/8)  \max( ||x||, |\eta| ) $. Now, fix $i$ such that 
$M=M_i$. 
By (\ref{eq:cat}) and the properties of $z_i$ we deduce that
\be \label{eq:fin}
{\cal S}  \geq
M   - 
\sum_{j \neq i } \Big( || K_j x +  \eta_j a_j||_i \Big) -
\frac{C_1}{8}||x||
\ee

We need to estimate $ || b_j||_i$ for 
$j \neq i$, where $b_j=K_j x + \eta a_j =  K_j z_j +  \eta_j a_j $.
We remark that, by the definitions of
$a_i$ and of $K_i$,  $ - \ddot{b}_j + A b_j =0$ on
$\R \backslash {\cal J}_i $. Hence 
$$
\max (| b_j (t) | , |\dot{b}_j (t)|/\lambda_2 ) \leq
e^{-\lambda_2 d(t, {\cal J}_j)} ||b_j||_j.
$$
So, for $i \neq j$, 
$$
||b_j||_i \leq e^{- \lambda_2 ( (T- T_{C_1}) + \ov{d}/2 )}
e^{-\lambda_2 |i-j-1| (\ov{d}+ 2 T})  ||b_j||_j. 
$$
Now,
$$
||b_j||_j \leq ||z_j||_j + ||z_j - K_j z_j - \eta_j a_j||_j
\leq ||x|| + M,
$$
where  we have  used  again that $||z_j||_j = ||x||_j $.
Therefore 
\begin{eqnarray*}
\sum_{j\neq i} ||K_j x + \eta_j a_j ||_i 
& \leq   &  e^{- \lambda_2 ( (T- T_{C_1}) + \ov{d}/2 )}
\sum_{j \neq i} e^{-\lambda_2 |i-j-1| (\ov{d}+ 2 T ) }  
(M+ ||x||)
 \\
&\leq& e^{- \lambda_2 ( (T- T_{C_1}) + \ov{d}/2 )}
 (1- e^{-\lambda_2 (\ov{d} + 2 T)})^{-1} 2 (M + ||x||)  \\
& \leq &\frac{7}{3} e^{- \lambda_2 ( (T- T_{C_1}) + \ov{d}/2 )}
(M +||x||),
\end{eqnarray*}
because $\ov{d} \geq 2/ \lambda_2$.
Combining this latter estimate and (\ref{eq:fin}), we get
\begin{eqnarray*}
{\cal S} 
& \geq & \Big( 1- \frac{7}{3} e^{- \lambda_2 ( (T- T_{C_1}) + \ov{d}/2 )} \Big)
M - \frac{7}{3} e^{- \lambda_2 ( (T- T_{C_1}) + \ov{d}/2 )} ||x||
-\frac{C_1}{8} ||x||
\\
&\geq & \Big( C_1 \Big[ 
\frac{7}{8} \Big( 1-\frac{7}{3} e^{- \lambda_2 ( (T- T_{C_1}) 
+ \ov{d}/2 )} \Big) 
- \frac{1}{8}\Big]- \frac{7}{3} 
e^{- \lambda_2 ( (T- T_{C_1}) + \ov{d}/2 )} \Big) 
||(x, \eta)||,
\end{eqnarray*}
which clearly implies the result in the lemma
(we use that $C_1 \leq 1$).

There remains to justify (\ref{eq:estK_i}).
Firstly, since $T> T_{C_1}$, ${\cal R}_i  x (s) =
\ov{\cal R}_i x (s) $ for $s \in {\cal J}_i$, 
${\cal R}_i x (s) =0$ for $s \notin {\cal J}_i$. Now, by 
$(W1)$, $(P1)$ and the definition of $T_{C_1}$, 
$|\ov{\cal R}_i x (s) | \leq \lambda_2^2 (C_1/8) ||x||$ for 
$s \notin {\cal J}_i$. Hence 
$|\ov{\cal R}_i x - {\cal R}_i x |_{\infty} \leq \lambda_2^2 
(C_1 / 8) ||x|| $, and by (\ref{eq:estL_A}) we get the 
first estimate in (\ref{eq:estK_i}).

Next we estimate $|{\cal R} x - \sum_{i=1}^k {\cal R}_i x|_{\infty}$.
On $[s_i, u_i]$ $Q_{\Theta} = Q^i$, hence 
$$
\forall s \in [s_i , u_i ] \quad
|{\cal R} x - \sum_{i=1}^k {\cal R}_i x| =
|\ov{\cal R}_i x - {\cal R}_i x | \leq \lambda_2^2 (C_1 /8) ||x||.
$$ 
So we have just to estimate $|{\cal R} x (s)|$ on $[s_i, u_{i+1}]$
for each $i \in \{ 1, \ldots , k-1 \}$. Note that 
on $[s_i , u_{i+1}]$  $Q_{\Theta}=
\gamma (Q^i (s_i) , Q^{i+1} (u_{i+1}), d_i)$.
A quite straighforward consequence of lemmas \ref{lem:homo} and
\ref{lem:link1} is that for $s \in [s_i, u_{i+1}]$,
$$
|Q_{\Theta} (s)| \leq 
\min \Big( \frac{7}{5} e^{\lambda_2 (s-u_{i+1})} + \frac{22}{21} 
e^{-\lambda_2 (s-s_i)} ,
\frac{7}{5} e^{-\lambda_2 (s-s_i)} + \frac{22}{21} 
e^{-\lambda_2 (s-u_{i+1})} \Big) r,
$$
$$
|\dot{Q}_{\Theta} (s)| \leq 
\min \Big( \frac{14}{5} e^{\lambda_2 (s-u_{i+1})} + \frac{22}{21} 
e^{-\lambda_2 (s-s_i)} ,
\frac{14}{5} e^{-\lambda_2 (s-s_i)} + \frac{22}{21} 
e^{-\lambda_2 (s-u_{i+1})} \Big) \lambda_1 r.
$$
Using that $d_i = u_{i+1}-s_i \geq 2/\lambda_2$ we can get
$$
\max \Big( |Q_{\Theta} (s)|, \frac{|\dot{Q}_{\Theta} (s)|}{\lambda_2}\Big)
\leq \frac{\lambda_1}{\lambda_2} r
\max \Big( \frac{14}{5} e^{-\lambda_2 d_i} + \frac{22}{21} , 
\Big( \frac{14}{5} + \frac{22}{21} \Big) 
e^{-\lambda_2 \frac{d_i}{2}} \Big)
\leq \frac{3\lambda_1}{2\lambda_2} r.
$$
Hence, by $(W1)$ and $(P1)$, for $s \in [s_i , u_{i+1}]$,
$|{\cal R} x (s)| \leq (3/2) \lambda_1 \lambda_2  \Lambda r ||x||$,
and by $(H4)$ 
$|{\cal R }x - \sum_{i=1}^k {\cal R}_i x |_{\infty}  \leq 
\lambda_2^2 (C_1/8) ||x||$. By (\ref{eq:estL_A}) this
implies the second estimate in (\ref{eq:estK_i}).
\end{pf}

\begin{lemma}{\label{lem:grad}} 
For $0 < r < r_1 $ and for all 
$\Theta = (\th_1, \ldots , \th_k ) \in \R^k $ with 
$\ov{d} >   2 / \lambda_2  $ there results that
$$
||S (Q_{\Th})|| \leq \frac{9}{2 \lambda_2} \max_{j=1,2} 
\Big\{  \lambda_j S_j e^{- \lambda_j \ov{d}}  \Big\} + 
7  r^2 \Lambda  e^{- \lambda_2 \ov{d}}.
$$
where $S_j = \max \{ | \ov{\al}_j | , | \ov{\bt}_j |, 
| \wtilde{\al}_j | , | \wtilde{\bt}_j | \}.$  
\end{lemma}

\begin{pf}
By construction  $Q_{\Th}$ solves  system (\ref{maineq}), 
except at  the times 
$  s_1, u_2,  \ldots  s_{k-1}, u_k  $  where 
$\dot{Q}_{\Th}$ is discontinuous and has a jump denoted by 
$\Delta \dot{Q}_{\Th}( u_i )$, $\Delta \dot{Q}_{\Th}( s_i )$. 
So we have for all $y \in Y \cap E$ 
$$
( S ( Q_{\Th}),y)= df(Q_{\Th})[y]=  \sum_{i=1}^{k-1} y( s_i  ) 
\cdot \Delta \dot{Q}_{\Th}( s_i ) + 
y( u_{i+1} ) \cdot \Delta \dot{Q}_{\Th}( u_{i+1} ).  
$$
Hence  $S(Q_{\Th})$ is the element of $E \cap Y$ defined by 
$$ 
 S(Q_{\Th}) =   
\sum_{i=1}^{k-1} g ( \cdot - s_i ) \Delta \dot{Q}_{\Th}( s_i )   +  
g ( \cdot - u_{i+1} ) \Delta \dot{Q}_{\Th}( u_{i+1} ),     
$$
where 
$- {\ddot g} + A g = \delta \ {Id}$, $\lim_{|t| \to \infty} g(t)=0$,
 $Id$ being the
identity  $2 \times 2$  matrix.
We have
$$
g (t)  = \frac{1}{2} (\sqrt{A})^{-1} e^{-\sqrt{A}|t|} \quad 
{\rm and } \quad \dot{g}(t) = - \frac{1}{ 2 } {\rm sign} (t) e^{-\sqrt{A}|t|}.
$$
We set $ \Delta = 
\max \{ |\Delta \dot{Q}_{\Th}( u_i )|, |\Delta \dot{Q}_{\Th}( s_i )| \ ;
\ 1 \leq i \leq k \}.$
We have to estimate $||S(Q_{\Th})||$.
Assume for example that $t \in [s_{n-1},u_n]$. We have
$$
|S(Q_{\Th})(t)| \leq \frac{\Delta}{2 \lambda_2} 
\Big ( 
\sum_{i=1}^{n-1} ( e^{-\lambda_2 | t - u_i| } + 
e^{-\lambda_2 | t - s_i| }) +  
\sum_{i=n}^k ( e^{-\lambda_2 | t - u_i| } 
+ e^{-\lambda_2 | t - s_i| }) \Big ).  
$$ 
For all $t\in [s_{n-1},u_n]$ we have
$$
\displaystyle
\frac{1}{2 \lambda_2 } \sum_{i=1}^{n-1} 
\Big( e^{-\lambda_2 ( t - u_i)} + e^{-\lambda_2 ( t - s_i )} \Big) 
\leq  \frac{1 + \exp{-( 2 \lambda_2 T) } }{2 \lambda_2 } 
\frac{\exp{-(t - s_{n-1})\lambda_2 }}{1-\exp{ -\lambda_2 ( \ov{d} + 2T)}}
 \leq  \frac{1}{ \lambda_2 }
\frac{\exp{-(t - s_{n-1})\lambda_2 }}{1-\exp{ -(\lambda_2 \ov{d} )}}.
$$
In the same way we get
$$
\sum_{i=n}^{k} \Big( e^{\lambda_2 ( t - u_i) } +
e^{\lambda_2 ( t - s_i )} \Big)  \leq  
\frac{1 + \exp{-( 2 \lambda_2 T ) } }{ 2 \lambda_2 }   
\frac{\exp{(t-u_n)} \lambda_2}{ 1 - \exp{- \lambda_2 ( \ov{d} + 2T  )}} 
\leq 
\frac{1}{  \lambda_2 } 
\frac{\exp{ (t-u_n) \lambda_2} }{ 1 - \exp{ - (\lambda_2  \ov{d}) } }.
$$
Hence for $t \in [s_{n-1}, u_n] $, since
$ \ov{d} > 2 / \lambda_2 $ we get
$ |S(Q_{\Th})(t)| \leq 
(\Delta (1 + \exp{-2} )) / (\lambda_2  (1 - \exp{-2} ))$. 
In the same way we can see that 
$|d/dt (S(Q_{\Th}))(t)| \leq 
(\Delta (1 + \exp{-2} )) / (\lambda_2  (1 - \exp{-2} ))4$. The case
$t \in [u_n , s_n]$ yields the same estimates. Hence
$$
|| S(Q_{\Th})|| \leq   
\frac{ \Delta (1 + e^{-2}) }{ \lambda_2 (1 - e^{-2}) }.
$$
We now estimate  $\Delta$.
By lemma \ref{lem:link1}  we have that 
$|\Delta - \Delta_L| \leq 5 \lambda_2 r^2 \Lambda e^{-\lambda_2 {\ov d}}. $
Hence  
\be\label{eq:final1}
|| S(Q_{\Th})|| \leq   
\frac{ (\Delta_L + (5 \lambda_2 r^2 
\Lambda e^{- \lambda_2 {\ov d}})) }{ \lambda_2} 
\frac{ 1 + e^{-2} }{ 1 - e^{-2} }.
\ee
Since $ u_n  -  s_{n-1}  > \ov{d} , \ \forall n $ 
by (\ref{eq:delta}) we have that:
\be\label{eq:f}
\Delta_L 
\leq  
2 \sqrt{2} \frac{1 + e^{-2} }{ 1 - e^{-4} }
\max_j \{  \lambda_{j} S_j e^{- \lambda_j \ov{d} } \}. 
\ee
By (\ref{eq:f}) and (\ref{eq:final1} )
we deduce the estimate of the lemma.

Now, since $\ov{d} > 2/\lambda_2$, $\lambda_1 e^{-\lambda_1 \ov{d}}
\leq \lambda_2 e^{-\lambda_2 \ov{d}}$, so, if (H4) is satisfied,
 using that
 $ S_1, S_2 \leq r  $, and $r \Lambda \leq 1 / 20 $
 we also have 
\be\label{eq:ult}
||S( Q_{\Th} )|| \leq 
 \frac{9}{2} r e^{- \lambda_2 \ov{d} } + 
7  r^2 \Lambda  e^{- \lambda_2 \ov{d} } \leq
5 
r e^{- \lambda_2 \ov{d} }.
\ee\end{pf}

In the next ``shadowing type''
lemma we repeat the  arguments of \cite{BB}
based on the contraction-mapping theorem 
in order to build the immersions ${\cal I}_k$. 

\begin{lemma}\label{lem:const}
Define 
$
D_2 = \frac{1}{\lambda_2} \ln 
\left( \frac{40 \ov{\Lambda} r}{C_1^2 } \right) $ and assume that (H4)
is satisfied. 
Then  $\forall \Theta=(\th_1 , \ldots , \th_k ) 
\in \R^k   $ with 
$ \ov{d} > \max \{ D_1, D_2 \}  $
there is a $C^0 $ function of $\Th$, $\Th \to  w( \Th )$ with 
$w(\Th)  \in Y $ such that:
\begin{itemize}
\item
$  S( Q_{\Th}  + w( \Th) )  \in  {\rm span} \{  a_1, \ldots , a_k  \} := 
\lla a_1 , \ldots , a_k \rn  $;
\item
$ w ( \Th )  \in 
\lla a_1 , \ldots , a_k \rn^{\bot}$.
\end{itemize}

Moreover ${\cal I}_k : 
\Teta \to Q_{\Teta} + w ( \Teta )$ is a $C^1$ function and
$$
|| w (\Th) ||  
\leq  \frac{18}{\lambda_2 C_1}     
\max_j \{ \lambda_j S_j \exp (- \lambda_j \ov{d}) \} +
\frac{28}{ C_1}  r^2 \Lambda  \exp (- \lambda_2 \ov{d}).
$$
In the sequel the function $w(\Th)$
 will be denoted also by  $ w_{ \Theta} $.   
\end{lemma}

\begin{pf}
This proof will follow closely the one given in \cite{BB},
see lemmas 3 and 13.
Therefore we shall be brief.
We shall use the following abbreviation: 
$$
F( \Th , w) = 
\frac{ \partial H}{ \partial 
(h, \mu )}_{|(\Th , Q_{\Th}+w)} \ \in \ L( Y \times \R^k ).
$$
By lemma \ref{lem:deriv} we know that 
$ \forall \Theta = (\th_1, \ldots , \th_k )  
\ |  \ \ov{d} > D_1, 
\ || F^{-1} ( \Th , 0) || \leq 2 /  C_1 .$

Let $B_{\delta} \subset Y \times \R^k $
be the ball in $ Y \times \R^k $ of center $0$ and radius $\delta$:
$B_{\delta} = \{( w , \mu_1 , \ldots ,  \mu_k )$  such that   
$ \max(||w|| , | \mu_1 | , \ldots , |\mu_k |) \leq \delta \}$.
We have to find $(w,\mu ) $ such that
$H (\th_1 , \ldots , \th_k ,Q_{\Th}+w, \mu_1 , \ldots ,  \mu_k ) = 0.$
This last equation is equivalent to  
$ {\cal D} ( w , \mu )  = (w , \mu) $ where: 
$$
{\cal D}(w, \mu ) =
 - F^{-1} ( \Th , 0) H(\Th, Q_{\Th}, 0) -
F^{-1} ( \Th , 0) 
\Big( H(\Th , Q_{\Th} +  w, \mu ) - 
H( \Th , Q_{\Th}, 0) - F ( \Th ,0) [w, \mu] \Big) .
$$
We will find  $ \delta  > 0 $ such that if 
$ \ov{d} > \max \{ D_1 , D_2  \}  $ then 
$$
{\bf (i)} \
\ov{{\cal D} (B_{\delta}) } \subset B_{\delta }, \qquad 
{\bf (ii)} \
{\cal D} \ {\rm is \  a  \ contraction  \ on} \ \ov{B_{\delta}}.
$$
It is easy to see that by  $(W1), (P1)$ and  (\ref{eq:estL_A}) 
$ || F( \Th , w) - F ( \Th , 0 ) || <  \ov{\Lambda} || w ||.$
 As in \cite{BB} we can derive that
$\forall (w, \mu ) \in B_{\delta}$ :
$$
||{\cal D}(w, \mu )|| \leq  \frac{2}{C_1} ||S( Q_{\Th} ) || + 
\frac{\ov{\Lambda}}{ C_1 }   ||(w, \mu)||^2  
$$
Then in order to get (i) we have to solve: 
\be\label{eq:dis}
\frac{2}{C_1} ||S( Q_{\Th} ) || + 
\frac{\ov{\Lambda}}{ C_1 }   \delta^2 < \delta  
\ee
A straigforward computation shows that if 
$
||S( Q_{\Th} ) || < C_1^2  / 8 \ov{\Lambda} 
$
then (\ref{eq:dis}) is satisfied for
$$
\delta \in 
\Big( C_1 
 \frac{  1 - \sqrt{ 1 - 8 \ov{\Lambda} ||S( Q_{\Th}) ||  / 
C_1^2 }}{ 2 \ov{\Lambda} }, 
C_1   \frac{ 1 + \sqrt{ 1 - 8 \ov{\Lambda} ||S( Q_{\Th}) ||  /
 C_1^2 }}{ 2 \ov{\Lambda}} \Big). 
$$ 
We now prove that also (ii) is satisfied if: 
$$
C_1 
 \frac{ 1 - \sqrt{ 1 - 8 \ov{\Lambda} 
||S ( Q_{\Th} )||  / C_1^2 }}{ 2 \ov{\Lambda}} 
< \delta < \frac {C_1}{2 \ov{\Lambda}}.
$$
Indeed, by $(W1)$, $(P1)$ and (\ref{eq:estL_A}), we have $ \forall (w, \mu ),
(w', \mu' ) \in B_{\delta}$:  
$$
 || {\cal D}( w, \mu ) - {\cal D}(w', \mu' ) || 
 \leq \frac{2}{C_1 }  \ov{\Lambda} \delta 
||(w, \mu )- (w', \mu' ) ||
$$
which implies our claim.
Then in order to apply the contraction-mapping theorem take
$\ov{d} > D_2 $ so that
$
|| S (Q_{\Th})|| <   C_1^2 / 8 \ov{\Lambda}.   
$
Taking into account (\ref{eq:ult}) 
the last inequality is satisfied if $\ov{d} > D_3 $ as
defined in the lemma.
Moreover, by the previous considerations, 
$$
||w_{\Th} || \leq C_1 
 \frac{  1 - \sqrt{ 1 - 8 \ov{\Lambda} ||S( Q_{\Th}) ||  /  
C_1^2 }}{ 2 \ov{\Lambda} }  \leq  \frac{4}{C_1} ||S (Q_{\Th})||,
$$
which, by lemma \ref{lem:grad} implies the  last estimate of the lemma.
The fact that $Q_{\Teta} + w(\Teta) $ is a $C^1$ function of 
$\Teta$ is a  consequence of the Implicit function theorem
applied to $H$.
\end{pf}

We define also 
$D_3 = \frac{1}{\lambda_2} \ln 
\left( \frac{40}{C_1 } \right) $ 
so  that for all $ \Theta \in \R^k $ with $ \ov{d} > \max \{ D_1, 
D_2, D_3 \} $ we have $||w(\Th)|| < r/2 $.

We define for $ \ov{d} > \ov{D}=\max \{ D_1 , D_2 , D_3  \} $  
the immersions ${\cal  I}_k $ 
$$ {\cal  I}_k  : \ M_k = \{ \Th \in \R^k \ | \ 
\th_{i+1} - \th_i > \ov{D} \} \to  
{\cal I}_k (\Th ) =   
Q_{\Th}  +  w ( \Th ) .
$$

By lemma \ref{lem:const} we can prove that:

\begin{lemma}\label{lem:nat}
If $\ov{d} > \ov{D}$ and (H4) is satisfied then
$ {\cal I}_k$ is a natural constraint for $f$.
\end{lemma}

The proof is in the appendix.

\subsection{Critical points of $f \circ {\cal I}_k$ and 
proof of theorem $1$}\label{sec:cri}

We are led, in order to find $k$-bumps homoclinics,
to look for critical points of the function 
$ f \circ {\cal I}_k ( \Th ) = f( Q_{\Th} + w ( \Th )) $.
Note that, since (\ref{maineq}) is autonomous,  
$ f( Q_{\Th} + w ( \Th )) $ depends only on  $d_1, \ldots, d_{k-1}$.
Let us define
$$
g( d_1 , \ldots , d_{k-1} ) = f( Q_{\Th} + w ( \Th )).  
$$
By lemma \ref{lem:nat} a  zero of the function 
$G : {\R}^{k-1} \to {\R}^{k-1}$ defined by
$$
G(d_1, \ldots , d_{k-1} ) = 
\Big( \frac{\partial g (d) }{\partial d_1} , \ldots ,
\frac{\partial g (d) }{\partial d_{k-1}} \Big)
$$
gives rise to an homoclinic solution of (\ref{maineq}). 
We will find a zero of $G$ by means of degree theory
showing in the proof of theorem \ref{thm:main2} that  
hypotheses $(H1-4)$ imply 
$$
|\deg (G ,U , 0) | =1
$$
where
$$
U = {\prod}_{i=1}^{k-1} ( D , J ) \subset \R^{k-1},
$$ 
and $J > D $ are some real numbers estimated in the proof of lemma
\ref{lem:link}.

We need  some preliminary lemmas. 
The next  one is proved in the appendix.

\begin{lemma}\label{lem:en}
For $d > 2 / \lambda_2 $ and $0 < r < r_1$ 
consider the solution $q_d$ given by lemma \ref{lem:link1} and 
the function 
$e(\bt, \al, d) $ which is the value of the action on $q_d$.
There results that 
$$
\frac{\partial e }{\partial d} ( \bt, \al , d ) =
- {\cal E}(d),
$$
where 
$ {\cal E} (d)= 
( {\dot q}^2_d (t) -
A q_d (t) \cdot q_d (t) )/2 + W (q_d) $
is the energy of the orbit $q_d$.
\end{lemma}

\begin{lemma}\label{lem:funct}
For all $(d_1, \ldots , d_{k-1} ) \in {\R}^{k-1} $ with 
$\ov{d} > \ov{D} $, we have:
\be\label{eq:deriv1}
\frac{\partial}{ \partial d_i } g(d_1, \ldots , d_{k-1} )=  
\left(\frac{\partial}{ \partial d } e \right) 
\Big( Q^i (s_i)+w_{\Th} (s_i),
Q^{i+1} (u_{i+1}) + w_{\Th} (u_{i+1}), d_i \Big).
\ee
\end{lemma}

\begin{pf}
We must compute
$$
\frac{ \partial g (d_1 , \ldots , d_{k-1} ) }{ \partial d_i } = 
\frac{ \partial f (Q_{\Teta} + w(\Teta) ) }{ \partial d_i }.
$$
For this purpose we consider the function of the real variable $ \tau $: 
$$
\sigma ( \tau ) = g (d_1 , \ldots , d_i + \tau , \ldots , d_{k-1}) = 
f( Q_{(\teta_1, \ldots , \teta_i, \teta_{i+1}+\tau , \ldots , 
\teta_k + \tau ) } +  w (\teta_1, 
\ldots , \teta_i, \teta_{i+1}+\tau , \ldots , 
\teta_k + \tau ))$$
and we compute $\sigma' (0)$. 
For simplicity of notation we set  $q^0 = Q_{\Teta} + w (\Teta)$. 
Let $q^{\tau} \in Y \cap E$ be defined as follows :
$$
\begin{array}{rcl}
q^{\tau} = 
\end{array}
\left\{
\begin{array}{rcl}
& & 
q^0 (s) \quad {\rm on } \quad (-\infty , s_i ]  \\
& & 
\gamma (q^0 (s_i ), q^0 (u_{i+1}) ,
d_i + \tau ) ( \cdot - s_i ) \ {\rm on } \ [s_i , u_{i+1} + \tau ] \\
& & 
q^0 ( \cdot - \tau ) \quad {\rm on } \quad [u_{i+1} + \tau , + \infty )
\end{array}\right.
$$
Note that our notations are coherent ( {\it i.e.} 
$q^{\tau} = q^0 $ when $\tau = 0$). 
We will use the notation  $a_i = a_{{\th}_i}$, not
distinguishing for simplicity between $a= \ov{a}$ and 
$a=\wtilde{a}$.  
\\[2mm]
Since $ q^0 - Q_{\Th} = w(\Teta) \in 
\lla a_{\teta_1} , \ldots , a_{\teta_k} \rn^{\bot}$ 
and since,  by the definition of $a_{\th_i}$,  
$ (x , a_{\th_i} ) = \alpha \int_{ J_i} x(t) \dot{Q}^i (t) dt  $ we see that
$$
q^{\tau} -  
Q_{ (\th_1 , \ldots , \th_i , \th_{i+1}+ \tau , \ldots ,\th_k + \tau )} 
\in \lla a_{\teta_1}, \ldots , a_{\teta_{i}},
a_{\teta_{i+1}+ \tau }, \ldots , a_{\teta_k + \tau } \rn^{\bot},
$$
and we can write
$ q^{\tau} - 
(Q_{(\teta_1, \ldots , \teta_i, \teta_{i+1} + \tau , \ldots , 
\teta_k + \tau ) } + w ( \th_1 , \ldots , \th_i, \th_i + \tau , \ldots ,
\th_k + \tau ))  =  \wtilde{w} (\tau),$
where $\wtilde{w} (0)=0$ and
$ \wtilde{w}(\tau) \in \lla a_{\teta_1}, \ldots , a_{\teta_{i}},
a_{\teta_{i+1}+ \tau }, \ldots , a_{\teta_k+ \tau}  \rn^{\bot}$.
Since $ \wtilde{w} (0) = 0 $ we have that 
$$ 
0 = (\partial / \partial {\tau })_{\tau=0}
  \ ( \wtilde{w}(\tau) , a_{\th_i + \tau} ) = 
( \partial \wtilde{w}(\tau ) / \partial {\tau} , 
a_{\th_i + \tau})_{\tau=0} + 
(\wtilde{w}(\tau) , \partial_{\tau} \ a_{\th_i + \tau})_{\tau=0} 
= (\partial \wtilde{w}(\tau ) / \partial {\tau} , a_{\th_i + \tau} 
)_{| \tau =0 }.
$$
This means that:
$$
\frac{ \partial \wtilde{w} ( \tau )}{  \partial {\tau}}_{| \tau=0 }  
\in \lla a_{\teta_1}, \ldots , a_{\teta_{i}},
a_{\teta_{i+1} }, \ldots , a_{\teta_k} \rn^{\bot}.
$$
Now we can prove (\ref{eq:deriv1}). Indeed, since 
$ (S (q^0), x)=0$ for all
$x \in \lla a_{\teta_1}, \ldots , 
 a_{\teta_k}   \rn^{\bot}$,
$$
\sigma'(0) =  \Big( S (q^0),  \frac{\partial ( q^{\tau} + 
\wtilde{w}(\tau))}{\partial \tau} \Big)_{| \tau = 0 } 
=  \left( S(q^0), \frac{\partial q^{\tau} }{\partial \tau}_{ | \tau = 0} 
\right) =  \frac{\partial f (q^{\tau} )}{\partial \tau}_{|\tau=0}.
$$
By the definition of $q^{\tau}$ we have that  
$$
\frac{\partial f (q^{\tau} )}{\partial \tau}_{|\tau=0} = 
\left( \frac{\partial }{ \partial d } e \right) 
 (q^0 (s_i)  , q^0 (u_{i+1}) ,  d_i),
$$
which yields (\ref{eq:deriv1}).
\end{pf}

\begin{lemma}\label{lem:ener}
For all $|\alpha|, |\beta| <  r' \leq r_1$, $ d > 2 / \lambda_2 $ 
there results that: 
\be\label{eq:deriv2} 
\left| 
\left( \frac{\partial}{ \partial d } e \right) 
( \bt , \al , d) -
\sum_{j=1}^2 \frac{\lambda_j^2 }{ (\sinh( \lambda_j d ))^2 } 
\left( \al_j \bt_j \cosh{(\lambda_j d) } - 
\frac{ ( \al_j^2 + \bt_j^2 ) }{2} \right) \right| \leq
\frac{15}{2} \lambda_1 \lambda_2 \Lambda r'^3 e^{-\lambda_2  d}.  
\ee
\end{lemma}

\begin{pf}
If  $0 < r' <r_1 $, by lemmas \ref{lem:homo} and \ref{lem:link1} 
the solutions $q_h^+$  and $q_d$ are defined. We will call 
$q_h = q_h^+$. 
Since ${\cal E}(d)$ is a constant of the  motion and
the homoclinic $q_h$ has zero energy we have
$$
{\cal E}(d) =
\frac{1}{2} (|{\dot q}_d(0)|^2  - A \bt \cdot \bt) + W( \bt ) \qquad
{\rm and} \qquad
0 = \frac{1}{2} (|{\dot q}_h(0)|^2  - A \bt \cdot \bt) + W( \bt ).
$$
Hence,  we obtain by substraction
\be\label{eq:en}
{\cal E}(d) = \frac{1}{2} ( |{\dot q}_d(0)|^2  - |{\dot q}_h(0)|^2 ).
\ee
Similarily for the linear system we have 
\be\label{eq:l}
{\cal E}_L(d) = \frac{1}{2} ( |{\dot q}_{d,L} (0)|^2 -
 |{\dot q}_{h,L} (0)|^2). 
\ee
This last expression can be  computed and we get 
$$
{\cal E}_L (d) = \sum_{j=1}^2 \frac{\lambda_j^2 }{ (\sinh( \lambda_j d ))^2 } 
\left(\al_j \bt_j \cosh{(\lambda_j d)}-\frac{(\al_j^2 +\bt_j^2 ) }{2} \right).
$$
By (\ref{eq:en}) and (\ref{eq:l}) we have
\begin{eqnarray*}
| {\cal E}(d) - {\cal E}_L(d)| & = &  
\frac{1}{2} \Big| (|{\dot q}_d (0)|^2 - |{\dot q}_h (0)|^2) - 
(|{\dot q}_{d,L}(0)|^2 - |{\dot q}_{h,L} (0)|^2) \Big| \\
& \leq & \frac{1}{2}
\Big| \Big( {\dot q}_d (0) - {\dot q}_h (0)) - 
({\dot q}_{d,L} (0) - {\dot q}_{h,L} (0) \Big) \Big| 
 \Big( | {\dot q}_d (0)|  + |{\dot q}_h (0)| \Big) \\
& + & 
\frac{1}{2} \Big| {\dot q}_{d,L}(0) - {\dot q}_{h,L} (0) \Big| 
\Big( | {\dot q}_d (0) - {\dot q}_{d,L} (0) | +
|{\dot q}_h (0)- {\dot q}_{h,L} (0) | \Big). 
\end{eqnarray*}
By lemmas \ref{lem:homo} and \ref{lem:link1} we have that 
\be \label{eq:diff1}
\Big| ({\dot q}_d (0) - {\dot q}_h (0)) - 
({\dot q}_{d,L}(0) - {\dot q}_{h,L} (0)) \Big| \leq 
5 \lambda_2 r'^2 \Lambda e^{- \lambda_2 d } \quad {\rm and } \quad
|{\dot q}_h (0)- {\dot q}_{h,L}(0) | \leq \frac{2}{7} \lambda_2 r'^2 \Lambda.
\ee
For $d \lambda_2 > 2 $ we get from (\ref{eq:diff1}) 
 that:
\be\label{eq:st}
|{\dot q}_d (0)- {\dot q}_{d,L} (0)|  \leq   \lambda_2  r'^2 \Lambda.
\ee
From (\ref{eq:diff1}) since
$r'\Lambda \leq 1 / 6 $ and $ \lambda_2 < \lambda_1 $ we obtain: 
\be\label{eq:derqh}
|{\dot q}_h (0)|  \leq \frac{22}{21} \lambda_1 r'. 
\ee
The expression of $q_{d,L}$, $q_{h,L}$ in subsection \ref{sec:bound} leads to
\be\label{eq:dli}
|{\dot q}_{d,L} (0)- {\dot q}_{h,L} (0)|  \leq  \frac{7}{3} r' \lambda_2 
e^{-\lambda_2 d }.
\ee
By (\ref{eq:st}) and (\ref{eq:linear}) we get
\be\label{eq:la}
|{\dot q}_d (0)|   \leq  \frac{4}{3} \lambda_1 r'
\ee
and finally by (\ref{eq:diff1}), (\ref{eq:st}), 
(\ref{eq:derqh}), (\ref{eq:dli}) and (\ref{eq:la}) we get:  
$$ 
| {\cal E}(d) - {\cal E}_L(d) | 
\leq \frac{15}{2} \lambda_1 \lambda_2 \Lambda  r'^3 e^{- \lambda_2 d },
$$
which is  (\ref{eq:deriv2}).
\end{pf}

The next  lemma is the most important for the proof of 
theorem \ref{thm:main2}.

\begin{lemma} \label{lem:link}
Assume $(H1-4)$. There exist $ \ov{D} < D < J $  such that 
for all $d = (d_1 , \ldots , d_{k-1} ) \in U$, 
 for all $ i \in \{ 1, \ldots , k-1 \} $  :
\begin{itemize}
\item 
if $ \alpha_1^{i+1} \beta_1^i >0 $ and $\alpha_2^{i+1} \beta_2^i <0 $,
$$ 
\partial_d e \Big( (Q_{\Teta} + w (\Teta ))(s_i) ,
(Q_{\Teta} + w (\Teta )) (u_{i+1}) , J \Big) < 0 
$$
and 
$$
\partial_d e \Big( (Q_{\Teta } + w (\Teta ))(s_i) ,
(Q_{\Teta} + w (\Teta )) (u_{i+1}) , D \Big) > 0. 
$$
\item 
if $\alpha_1^{i+1} \beta_1^i <0 $ and $\alpha_2^{i+1} \beta_2^i >0 $,
$$
\partial_d e \Big( (Q_{\Teta} + w (\Teta ))(s_i) ,
(Q_{\Teta} + w (\Teta )) (u_{i+1}) , J \Big) > 0 
$$
and
$$
\partial_d e \Big( (Q_{\Teta} + w (\Teta ))(s_i) ,
(Q_{\Teta} + w (\Teta )) (u_{i+1}) , D \Big) < 0. 
$$
\end{itemize}
\end{lemma}

Estimates for $D$ and $J$ are given in the proof.

\begin{pf}
By $(H1)$ the following cases can arise:
$\alpha_1^{i+1} \beta_1^i >0 $ and $\alpha_2^{i+1} \beta_2^i <0 $
or
$\alpha_1^{i+1} \beta_1^i <0 $ and $\alpha_2^{i+1} \beta_2^i >0 $.
We first deal with the first case.
Set 
$$
(\hat{\bt}_1^i, \hat{\bt}_2^i ) =
(Q_{\Teta} + w (\Teta ))(s_i)
= (\bt_1^i, \bt_2^i) + w(\Teta) (s_i) 
$$
and
$$
(\hat{\al}_1^{i+1}, \hat{\al}_2^{i+1} ) =
 (Q_{\Teta} + w (\Teta ) ) (u_{i+1} ) 
= (\al_1^{i+1} , \al_2^{i+1} ) + w(\Teta)(u_{i+1}).
$$
We will choose $D$ large enough such that 
\be \label{eq:comp}
\lambda_1 S_1 e^{-\lambda_1 D} \leq 
2 \lambda_2 S_2 e^{-\lambda_2 D}.
\ee
So, by lemma \ref{lem:const}, 
$$
|\hat{\alpha}- \alpha|, |\hat{\beta}-\beta| \leq 
||w(\Teta)|| \leq
\frac{36 S_2 + 28 r^2 \Lambda}{C_1} e^{-\lambda_2 D}.
$$
One of our conditions will be
\be \label{eq:condep} 
\Big| \frac{36 S_2 + 28 r^2 \Lambda }{C_1} e^{-\lambda_2 D} \Big|
\leq \ep {\cal M},
\ee
for some $\ep \in (0,1/8)$ which will be chosen later on.
Note that  $D \geq D_3$ then holds. This implies also
\be \label{eq:apep}
\hat{\alpha} \in [(1-\ep)\alpha , (1+\ep)\alpha ], \
|\hat{\alpha}| \leq (1+\ep)r  \quad 
; \quad \hat{\beta} \in [(1-\ep)\beta , (1+\ep)\beta ],
\ |\hat{\beta}| \leq (1+\ep)r.
\ee
Then there results that
$ \hat{\alpha}_1^{i+1} \hat{\beta}_1^{i} > 0 
\quad  {\rm and} \quad
 \hat{\alpha}_2^{i+1} \hat{\beta}_2^{i} < 0.$ 
Now, by lemma \ref{lem:ener}, we can write 
$$
\frac{\partial e}{\partial d}( \hat{\beta}^i , \hat{\al}^{i+1}, d_i )  =
-U_2 \frac{\cosh (\lambda_2 d_i)}{(\sinh (\lambda_2 d_i))^2}
+U_1 \frac{\cosh (\lambda_1 d_i)}{(\sinh (\lambda_1 d_i))^2},
$$
where
$$
U_2 \leq \lambda_2^2 \Big ( 
| \hat{\alpha}_2^{i+1} || \hat{\beta}_2^i  | + 
\frac{ (\hat{\alpha}_2^{i+1})^2 +  (\hat{\beta}_2^i )^2}{ 2 
 \cosh (\lambda_2 d_i) } \Big ) +
\frac{15}{4} \lambda_1 \lambda_2  (1+ \ep)^3 r^3 \Lambda 
$$
and
$$
U_1 = \lambda_1^2 \Big ( 
| \hat{\alpha}_1^{i+1} || \hat{\beta}_1^{i+1} | - 
\frac{ ( \hat{\alpha}_1^{i+1} )^2
+ ( \hat{\beta}_1^i )^2 }{2 \cosh (\lambda_1 d_i) } \Big ).
$$
By (\ref{eq:apep}) and (\ref{eq:condep}) we get
$$
U_2 \leq \lambda_2^2  (1+\ep)^2 \Big (
|\alpha_2^{i+1}||\beta_2^i| 
+ \Big( (\alpha_2^{i+1})^2
+(\beta_2^{i+1})^2 \Big) \frac{C_1 \ep {\cal M}}{36 S_2}\Big ) +
\frac{15}{4} \lambda_1 \lambda_2 (1+\ep)^3 r^3 \Lambda ,
$$
hence, since $C_1 \leq 1$, ${\cal M} \leq \min \{
|\alpha_2^{i+1}| , |\beta_2^i | \}$ and $S_2 \geq \max \{
|\alpha_2^{i+1}| , |\beta_2^i | \}$, 
\be \label{eq:u2}
U_2 \leq \lambda_2^2  (1+\ep)^2 (1 + \frac{ \ep}{9}) 
|\alpha_2^{i+1}||\beta_2^i| +  
\frac{15}{4} \lambda_1 \lambda_2 (1+\ep)^3r^3 \Lambda. 
\ee
On the other hand,  (\ref{eq:condep}) and (\ref{eq:comp}) imply that
$$
\frac{18 \lambda_1}{C_1 \lambda_2} S_1 e^{-\lambda_1 D}
\leq \ep {\cal M}.
$$
We derive readily that
$$
\frac{ (\hat{\alpha}_1^{i+1})^2
+ (\hat{\beta}_1^i)^2}{2 \cosh (\lambda_1 d_i)} <
\frac{C_1 \lambda_2}{9 \lambda_1}
(1+\ep)^2 |\alpha_1^{i+1}||\beta_1^i|,
$$
hence, by (\ref{eq:apep}),
\be \label{eq:u1}
U_1 > \lambda_1^2 (1-\ep)^2 
\Big( 1-\frac{C_1 \lambda_2}{9 \lambda_1}(\frac{1+\ep}{1-\ep})^2 \Big) 
|\alpha_1^{i+1}||\beta_1^i|.
\ee
From (\ref{eq:u2}),  (\ref{eq:u1}) and the fact that 
$\ep < 1/8$  we get
\be \label{eq:u1u2}
\frac{U_2}{U_1} < \frac{\lambda_2^2 (1+\ep)^3 }{ \lambda_1^2
(1-\ep )^2 (1-\ep/5) } 
\Big( 
\frac{|\alpha_2^{i+1}||\beta_2^i|+ (15 \lambda_1 /
4 \lambda_2) \Lambda r^3 }{|\alpha_1^{i+1}||\beta_1^i |} \Big)
\ee
We shall take 
\be \label{eq:defd}
D= \min_i  \frac{1}{\lambda_1 -\lambda_2}
\ln \left(\frac{ \lambda_1^2
(1-\ep )^2 (1- \ep /5) }{\lambda_2^2 (1+\ep)^3 }
\Big(  \frac{ | \alpha_1^{i+1} | | \beta_1^i| }{ 
|\alpha_2^{i+1}||\beta_2^i|+ (15 \lambda_1 / 4 \lambda_2) \Lambda r^3 }
 \Big) \right). 
\ee
It is easy to see that, if $d_i =D$ and
$d_j \geq D$ for all $j$ then, by (\ref{eq:u1u2}),
$$
\frac{\partial e}{\partial d} ( \hat{\beta}^i , \hat{\al}^{i+1}, d_i)  > 0 
$$
Now,  as $ d_i \to + \infty $,
$ \partial_d e( \hat{\beta}^i , \hat{\al}^{i+1} , d_i ) 
 \wtilde{=}  -2U_2 e^{-\lambda_2 d_i},  $
where
$$
U_2 \geq \lambda_2^2 (1-\ep)^2 (1- \ep/5) |\alpha_2^{i+1}||\beta_2^i| -
\frac{15}{4} \lambda_1 \lambda_2 (1+\ep)^3  r^3 \Lambda
$$
(This estimate could be improved in the first case but we want 
to be able to extend our arguing to the second case $\alpha_2^{i+1}
\beta_2^i >0$).
Hence, provided 
\be \label{eq:condr}
|\alpha_2^{i+1}||\beta_2^i| > 20 \lambda_1  r^3 \Lambda / \lambda_2   ,
\ee
we get
$\partial_d e ( \hat{\beta}^i , \hat{\alpha}^{i+1} , J ) < 0,$
for  $J$  large enough, more exactly for
$$
J > \max_i  \frac{1}{\lambda_1 -\lambda_2 }
\ln \left( \frac{\lambda_1^2}{\lambda_2^2}
\frac{|\alpha_1^{i+1}||\beta_1^i|}{|\alpha_2^{i+1}||\beta_2^i|- 
(20 \lambda_1 \Lambda r^3 / \lambda_2 ) }\right).
$$
Therefore we get the desired result provided 
 conditions (\ref{eq:comp}), (\ref{eq:condep}) and (\ref{eq:condr}) 
 are satisfied, with
$D$ defined by (\ref{eq:defd}). Now we must choose $\ep$ to make
 condition (\ref{eq:condep})  as weak as possible.
This condition reads
\be \label{condfin}
\frac{\lambda_2^2  }{ \lambda_1^2 }
 \frac{ |\alpha_2||\beta_2|+ (15 \lambda_1 / 4 \lambda_2) \Lambda r^3}{|\alpha_1||\beta_1|}
\leq \frac{(1-\ep)^2 (1-\ep/5)}{(1+\ep)^3} 
\ep^{(\lambda_1/\lambda_2) -1} 
\left( \frac{C_1 {\cal M} }{ 36 S_2 + 28 \Lambda r^3}
\right)^{(\lambda_1/\lambda_2) -1}.
\ee
Therefore  we get condition $(H3)$.
We get condition $(H2)$ so that lemmas \ref{lem:grad} 
and \ref{lem:const} are satisfied with
our choice of $D$. The first inequality in 
condition $(H4)$ is just (\ref{eq:condr}). 

There remains to check that $(H2-4)$ imply that 
(\ref{eq:comp}) holds true. First (by $(H2)$ and the definition 
of $D$ in (\ref{eq:defd}), $D \geq 2/\lambda_2$, hence
$\lambda_1 e^{-\lambda_1 D} \leq \lambda_2 e^{-\lambda_2 D}$, and if
$S_1 \leq 2 S_2$ then (\ref{eq:comp}) holds. So we shall
assume that $S_1 > 2 S_2$. Fix $i$ such that 
\be \label{defd2}
D=   \frac{1}{\lambda_1 -\lambda_2}
\ln \left(\frac{ \lambda_1^2
(1-\ep )^2 (1- \ep /7) }{\lambda_2^2 (1+\ep)^3 }
\Big(  \frac{ | \alpha_1^{i+1} | | \beta_1^i| }{ 
|\alpha_2^{i+1}||\beta_2^i|+ (15 \lambda_1 / 4 \lambda_2) \Lambda r^3 }
 \Big) \right).
\ee
Note that, by $(H4)$, $(\Lambda \lambda_1 r^3/\lambda_2) \leq
|\al_2^{i+1}||\beta_2^i|/20$. Combined with the fact that 
$\ep \in (0, 1/8)$, this readily implies that
$e^{(\lambda_1 - \lambda_2) D} \geq 2 \lambda_1^2 
|\alpha_1^{i+1}||\beta_1^i| (5 \lambda_2^2 |\alpha_2^{i+1}|
|\beta_2^i|)^{-1}$.

Now $S_2 \leq S_1/2 \leq r/2$, hence
$|\alpha_2^{i+1}| \leq r/2$ and $|\alpha_1^{i+1}| =
(r^2 - (\alpha_2^{i+1})^2 )^{1/2} \geq \sqrt{3}r/2 $;
$|\beta_1^i| \geq \sqrt{3}r/2$ as well. So we get
$$
\frac{\lambda_2 S_2}{\lambda_1 S_1} e^{(\lambda_1 - \lambda_2 ) D}
\geq \frac{\lambda_2 S_2}{\lambda_1 S_1} \frac{2}{5}
\frac{3 r^2}{|\alpha_2^{i+1}|
|\beta_2^i| 4} \geq \frac{3 \lambda_2 r}{10 \lambda_1 S_2}
\geq \frac{3}{5} \geq \frac{1}{2}, 
$$
since $S_2 \leq r/2$. So (\ref{eq:comp}) holds.

We have proved the lemma in the case where 
$\alpha_1^{i+1} \beta_1^i > 0$ and $ \alpha_2^{i+1} 
\beta_2^i < 0$. The second case can be dealt with in
a similar way (in fact the estimates are simpler in that 
 case).
\end{pf}

We now show, using the previous lemma how to prove theorem 
\ref{thm:main2}.
\\[2mm] 
\begin{pfn}{\sc of theorem \ } \ref{thm:main2} .
Let $J$ be given by
lemma \ref{lem:link}. 
By lemmas \ref{lem:funct} and \ref{lem:link}, 
$g$ has the following
property :
for all $i \in \{ 1, \ldots , k-1 \}$, we have either
\\[2mm]
{\it $(P_-)$ : For all $d= (d_1 , \ldots , d_{k-1} )
\in \ov{U}$,
$$
d_i =D \Rightarrow 
\frac{\partial g}{\partial d_i} (d) > 0
\quad {\rm and} \quad
d_i =J \Rightarrow 
\frac{\partial g}{\partial d_i} (d) < 0,
$$} 
or \\[2mm]
{\it $(P_+)$ :
For all $d= (d_1 , \ldots , d_{k-1} )
\in \ov{U}$, 
$$
d_i =D \Rightarrow 
\frac{\partial g}{\partial d_i} (d) < 0
\quad {\rm and} \quad
d_i =J \Rightarrow 
\frac{\partial g}{\partial d_i} (d) > 0.
$$ }
It can be readily seen that this 
property implies $|\deg (G ,U , 0) | =1.$
In fact define the function $\hat{G} : \ov{U} \to \R^{k-1}$ as follows :
$$
\hat{G}(d_1, \ldots , d_{k-1}) =
\Big( \ep_1 (d_1 - \frac{D+J}{2} ), \ldots,
 \ep_{k-1} (d_{k-1} - \frac{D+J}{2} ) \Big) ,
$$
where $\ep_i =1$ if $(P_+)$ is satisfied for the
index $i$, and $\ep_i = -1$ if 
$(P_-)$ is satisfied for the index $i$. 
Since the homotopy $G_t = (1-t) \hat{G} + t G$ for $t \in [0,1]$ 
is admissible there results that 
$ \deg (G,U,0)$  $=$ $ \deg (\hat{G},U, 0)=$ $ \pm 1$ 
and the existence of a critical point of $g$ in $U$ follows.
This critical point corresponds to a homoclinic, which, 
by (\ref{eq:condep}) and since $\ep \in (0, 1/8)$, 
 enjoys the properties given in theorem \ref{thm:main2}.
\end{pfn}

\section{Examples}\label{sec:exam}
The aim of this section is to show examples of 
Hamiltonian systems where  the hypotheses 
$(H1-4)$ can be checked.

\subsection{Almost equal eigenvalues}\label{sec:almost}
Consider the following system  
$$
-\ddot{q}+\psi (q){\cal J} \dot{q} +A_{\ep}q = \nabla W (q),
\leqno{(S_{\ep})}
$$
with $A_{\ep}=
\left( \begin{array}{cc}
(\lambda +\ep)^2 & 0 \\
0   & (\lambda -\ep )^2 
\end{array} \right)$ and $\lambda >0 $. We assume that $W, \psi$ satisfy
(W1), (P1).
We shall use the following assumptions: 
\begin{itemize}
\item
(A1) $(S_0)$ has two nondegenerate homoclinics 
$\ov{q}$ and $\wtilde{q}$. 
\end{itemize}
It can be shown that the limits as
$t\to +\infty$ and as $t\to -\infty$ of
$\ov{q}(t) / |\ov{q}(t)|$ (resp.
$\wtilde{q}(t) / |\wtilde{q}(t)|$) do exist. Call 
$(\cos \ov{\om}_s , \sin \ov{\om}_s)$ and
$
(\cos \ov{\om}_u,\sin \ov{\om}_u)$
(resp. $(\cos \wtilde{\om}_s , \sin \wtilde{\om}_s)$ and
$(\cos \wtilde{\om}_u,\sin \wtilde{\om}_u)$) these limits. The second
assumption is 
\begin{itemize}
\item (A2)
$ {\om}_u, {\om}_s
\neq n \pi/2 , \  n \in \Z\  $, 
$\ -1 < \tan  \om_u   \tan  \om_s  < 0$ 
and ( $ \cos \ov{\om}_u  \cos \wtilde{\om}_u < 0 $ or
$\cos \ov{\om}_s  \cos \wtilde{\om}_s < 0$).
\end{itemize} 
As an application of theorem \ref{thm:main2} we get
\begin{theorem}\label{thm:alm}
Assume that $(S_0)$  satisfies
assumptions $(A1)$ and $(A2)$. Then there is $\ep_1 >0$ such that,
for $0<|\ep| < \ep_1$, $(S_{\ep})$ has a rich family of 
homoclinics, which induces a chaotic behaviour at the 
zero energy level, according to thm. \ref{thm:main2}. 
Moreover there is $C>0$ such that 
$h_{top}^0>C\ep$, where $h_{top}^0$ denotes the topological
entropy at  the zero energy level.
\end{theorem} 

\begin{pf}
Call $\ov{q}_0$ and $\wtilde{q}_0$ the two nondegenerate 
homoclinics for $(S_0)$. By the Implicit function theorem,
it is easy to get, for $|\ep|$ small enough, the existence 
of two homoclinics $\ov{q}_{\ep}$ and $\wtilde{q}_{\ep}$ for system 
$(S_{\ep})$,
and to see that these two homoclinics are nondegenerate, 
of constant of nondegeneracy $C_{1,\ep}$ which tends to
$C_1$ as $\ep \to 0$. Moreover
\be \label{eq:eginf}
\lim_{\ep \to 0} \max \Big \{ |\ov{q}_{\ep} - \ov{q}|_{\infty},
|\wtilde{q}_{\ep} - \wtilde{q}|_{\infty} \Big \} =0.
\ee
There are $r_2 >0$ and $\ep_0 >0$ such that,  for $0< r < r_2 $
and $|\ep| < \ep_0$,
the trajectory of $\ov{q}_{\ep}$ (resp. 
$\wtilde{q}_{\ep}$) crosses the circle of radius
$r$ at two points only:  
$\ov{\alpha}^{\ep} (r)$ and $\ov{\beta}^{\ep}(r)$ 
(resp. $\wtilde{\alpha}^{\ep} (r)$ and $\wtilde{\beta}^{\ep}(r)$).
Moreover, 
$$
\lim_{r \to 0} \frac{\ov{\alpha}^0 (r)}{r }=
(\cos \ov{\om}_u , \sin \ov{\om}_u ) ; \quad
\lim_{r \to 0} \frac{\ov{\beta}^0 (r)}{r }=
(\cos \ov{\om}_s , \sin \ov{\om}_s )
$$
and for all $r \in (0,r_2 )$, from (\ref{eq:eginf}),
\be \label{eq:eglim}
\lim_{\ep \to 0} \ov{\alpha}^{\ep}(r)
=\ov{\alpha}^{0}(r) ; \quad
\lim_{\ep \to 0} \ov{\beta}^{\ep} (r)
=\ov{\beta}^{0} (r).
\ee
We have similar properties for $\wtilde{q}, \wtilde{q}_{\ep}$.
By $ (A2) $, there are  $0< r_3 < r_2$ and 
$\delta >0$ such that $r_3 < \min \{ C_1 / (24 \Lambda), \rho_0 /2 \}$
and
\be \label{eq:egal2}
\alpha_2^0 (r_3)^2  \geq 40 \frac{\lambda_1}{\lambda_2} \Lambda 
r_3^3 + \delta r_3^2 ; \quad
\beta_2^0 (r_3)^2  \geq  40 \frac{\lambda_1}{\lambda_2} \Lambda 
r_3^3 + \delta r_3^2 ,
\ee
\be \label{eq:egal1}
|\alpha_1^0 (r_3)| \geq  \delta r_3 ; \quad
|\beta_1^0 (r_3)| \geq  \delta r_3 
\ee
and 
\be \label{eq:egtan}
\frac{|\alpha_2^0 (r_3)||\beta_2^0 (r_3)| + 
(15 \Lambda r_3^3 / 4) }{|\alpha_1^0
(r_3)|
|\beta_1^0 (r_3)|} < 1-\delta,
\ee
where $ \alpha_i $ (resp. $ \beta_i$) may represent either 
$ \ov{\alpha}_i$ or $\wtilde{\alpha}_i$ (resp.
$ \ov{\beta}_i$ or $\wtilde{\beta}_i$).

Now the eigenvalues of the equilibrium $0$ for 
$(S_{\ep})$ are $\pm \lambda_{2,\ep} =
\pm ( \lambda -|\ep|) $ and $\pm \lambda_{1,\ep} = \pm
( \lambda +|\ep|)$, hence 
$\lim_{\ep \to 0} \lambda_{1,\ep}/\lambda_{2,\ep} =1$.

From (\ref{eq:eglim}),(\ref{eq:egal2}) and (\ref{eq:egal1}) it is easy
to see that all the second members in conditions $(H2-3)$ 
(associated to $(S_{\ep})$) tend to $1$
as $\ep \to 0$. Therefore ( taking $r =r_3$), 
there is $\ep_1 > 0$ such that  these conditions and condition (H4) 
are satisfied,  by (\ref{eq:egtan}),
(\ref{eq:eglim}), (\ref{eq:egal2}) and 
(\ref{eq:egal1}),    for $ 0<|\ep | < \ep_1$. 
By theorem \ref{thm:main2} there is 
chaos at the zero energy level for 
$0<|\ep|< \ep_1$. The estimate on the topological entropy follows 
by the results of  section \ref{sec:dyn} since the distance between 
two consecutive bumps is of order $ 1 / \ep $
(see also the proof of the relaxed theorem \ref{thm:rela}). 
\end{pf}
\subsection{Perturbation of an uncoupled system}\label{sub:uncoup}

Let us consider a perturbed system of the following form
\be
\begin{array}{rcl}
{\label{eq:pert}} 
- {\ddot q}_1 + \lambda_1^2 q_1 & = & W_1'(q_1) + 
\ep \psi (q)  \dot{q}_2  \\
- {\ddot q}_2 + \lambda_2^2 q_2 & = & W_2' ( q_2 )  - 
\ep \psi (q)  \dot{q}_1
\end{array}
\ee
with $( q_1 , q_2 ) \in \R^2 $. We assume that $W_i (0)=
W'_i (0)= W''_i(0) =0$ for $i=1,2$ and that $\psi (0) =0$.
(\ref{eq:pert}) can be written as
\be \label{eq:pert2}
-\ddot{q} +  \ep \psi (q) {\cal J} \dot{q} + Aq = \nabla W (q) ,
\ee
where $W(q_1 , q_2) = W_1 (q_1) + W_2 (q_2)$.

For $ \ep = 0 $ system (\ref{eq:pert}) splits into 
the direct product of two $1$-dimensional systems.

For the sake of simplicity we shall suppose that 
$W_1$ and $W_2$ are even, and that 
$$
\psi (q_1 , -q_2) = \psi (-q_1 , q_2) = \psi (q_1 , q_2 ). 
$$
As a consequence, if $q=(q_1 (t) , q_2 (t) ) $ is a 
homoclinic solution to
(\ref{eq:pert2}), then $(q_1 (-t) , -q_2 (-t)) $ and $-q (t) $ 
are homoclinic solutions as well.

Suppose  that:
\be\label{eq:oned}
- {\ddot q}_1 + \lambda_1^2 q_1  =  W_1'(q_1)
\ee
possesses an homoclinic $q_0$. Up to  a time translation,
we may assume that $q_0$ is even.  
Thus,   for $ \ep = 0 $, 
  $\ov{q} = ( q_0 , 0 )  $ and 
$ \wtilde{q}= ( - q_0 , 0 ) $ are two nondegenerate
(up to time translation) homoclinic solutions of (\ref{eq:pert}). 
We define
$$
\Gamma = \int_{-\infty}^{+\infty} - \psi (q_0 (s), 0) 
\dot{q}_0 (s) e^{\lambda_2 s} \; ds.
$$
Note that we have $|q_0 (t)| +  |\dot{q}_0 (t)| \leq C e^{-
\lambda_1 |t|} $ for some positive constant $C$, hence, by the properties 
of $\psi$, $|\psi (q_0 ,0) \dot{q}_0 (t)| \leq C' e^{-2\lambda_1 |t|}$,
so the integral $\Gamma$ is well defined.  
As an  application of theorem \ref{thm:main2} we get:
 
\begin{theorem} \label{th:pert}
If $\Gamma \neq 0$, then there is $\ep_0 > 0$ such that, 
for $\ep \in (-\ep_0,0 ) \cup (0, \ep_0)$, (\ref{eq:pert})
has a rich family of homoclinics and a chaotic behaviour at the 
zero energy level. 
\end{theorem} 

Before proving this theorem we introduce $h_2$, defined by
$$
- {\ddot h}_2 + \lambda_2^2 h_2  = 
- \psi (q_0 , 0 )  \dot{q_0} \qquad {\rm and } \qquad
\lim_{|t| \to +\infty } h_2 (t) = 0.
$$
Solving  this  equation  we find:
\be\label{eq:expre}
h_2 (t) =  \frac{1}{2 \lambda_2 }
\left( \int_{t}^{+\infty} f(s) \exp (- \lambda_2 s ) ds 
\right) e^{ \lambda_2 t } + 
\frac{1}{ 2 \lambda_2 } 
\left( \int_{-\infty}^{t} 
f(s) \exp ( \lambda_2 s ) ds \right) e^{ - \lambda_2 t } 
\ee
where $f(s)= \psi ( q_0 (s ) , 0 )  \dot{q}_0 ( s ) $ 
is an odd function.
It is easy to see that $h_2 (t) \sim \Gamma e^{-\lambda_2 t} / 2\lambda_2 $
 as $t \to +\infty$ and that 
$h_2 (t) \sim -\Gamma e^{\lambda_2 t} / 2 \lambda_2 $
as $t \to - \infty$. We have

\begin{lemma} \label{lem:pert}
There are $\ep_1 > 0$ and a non increasing  function $a(\ep)$ which
tends to $0$ as $\ep$ tends to $0$  such that, for all 
$\ep \in (-\ep_1 , \ep_1)$, (\ref{eq:pert}) has a homoclinic 
solution $\ov{q}_{\ep}= (\ov{q}_{1, \ep}, 
\ov{q}_{2, \ep})$  satisfying 
$|\ov{q}_{1, \ep} - q_0 | \leq a(\ep) \ep e^{-\lambda_2 |t|}$,
$|\ov{q}_{2, \ep} - \ep h_2 | \leq a(\ep) \ep  e^{-\lambda_2 |t|}$,
$\ov{q}_{1, \ep} (-t) = \ov{q}_{1, \ep} (t)$,
$\ov{q}_{2, \ep} (-t) = -\ov{q}_{2, \ep} (t)$.
\end{lemma}
\begin{pf}
This is a consequence of the Implicit function theorem. Define 
the Banach space
$$
X'= \{ q=(q_1,q_2) \in X \  |  \  
q_1 (-t) = q_1 (t) \  ;  \  q_2 (-t) = - q_2 (t) \}.
$$
$X'$ is endowed with norm  $ || \cdot ||$  defined by
$||q||_1 = \max ( |q e^{\lambda_2  |t|}|_{\infty} , 
|\dot{q} e^{\lambda_2 |t|}|_{\infty} / \lambda_2 )$.

Let $F(q) = L_A ( \nabla W(q))$ and
$G(q) = L_A ( \psi (q) {\cal J} \dot{q})$. It is easy to see, by the
properties of $\psi$ and $W$, that :
\begin{itemize}
\item $F$ and $G$ map $X'$ into itself;
\item $F,G : X' \to X'$ are smooth and 
$dF ((q_0 , 0)), dG ((q_0 , 0)) $ are compact linear operators.
\end{itemize}
 (\ref{eq:pert2}) is equivalent to 
\be \label{eq:pert3}
q= F(q) - \ep G(q)
\ee
Now, it is a standard fact that for $\ep = 0$ the linearization of
(\ref{eq:pert}) at $(q_0 , 0)$ has no other homoclinic
solution $ (v_1 , v_2 )$ that satisfies $ v_1 (-t) = v_1 (t) \ ,
\ v_2 (-t) = - v_2 (t) $ than $0$. Hence 
$ Ker (I - dF ( q_0 , 0)) = 0$ and, since 
$ dF ((q_0 , 0)) $ is compact, $ I - dF (( q_0 , 0))   $ is 
an isomorphism from $X' $ to $X'$. Therefore we may apply the
Implicit function theorem and we get, for all $\ep$ small enough in
modulus, a solution $\ov{q}_{\ep}$ of (\ref{eq:pert}). Moreover
$(d/d\ep)_{\ep=0} ( \ov{q}_{\ep}) =(0, h_2)$, which is the solution of the 
linear equation $ (I - dF (q_0, 0) )h = - G(q_0 , 0)$. From this 
the estimates of the lemma follow.  
\end{pf}

{\bf Proof of theorem \ref{th:pert}}: 
Without loss of generality we assume that $ \Gamma < 0$ 
we perform the proof for $\ep >0$. Then, by lemma
\ref{lem:pert} there are $\ep_2$ and $\hat{T}<0$ such that, for 
$0< \ep < \ep_2$, 
and $t\leq \hat{T}$, $q_{2, \ep} (t) > 0$. We shall
prove the following lemma:
\begin{lemma}
For all $\omega$ small enough there is $0< \ep ( \omega) < \ep_2$ 
such that, for all $\ep \in (0 , \ep (\omega))$,
there are $T_{\ep}\leq \hat{T}, r_{\ep}$ 
such that $|\ov{q}_{\ep} (T_{\ep}) |
= r_{\ep}$, $ |\ov{q}_{\ep} (t) | < r_{\ep}$ for $t< T_{\ep}$,
and $q_{2, \ep} (T_{\ep}) / q_{1, \ep } (T_{\ep}) = \tan \omega$.
In addition $\lim_{\ep \to 0} T_{\ep} = - \infty ,
\lim_{\ep \to 0} r_{\ep} = 0$. 
\end{lemma}

\begin{pf}
This is a  consequence of lemma \ref{lem:pert}. We may assume
that $\ep_1$ and $\omega$ are small enough such that
$a(\ep_1) / (|\Gamma| - a (\ep_1)) < \cot \om$. 
For $0< \ep <  \ep_1 $ and $t \leq \hat{T}$, define $f_{\ep} (t) $ by
$f_{\ep} (t) = q_{1, \ep}(t)  / q_{2, \ep} (t)$.
This is a continuous function on $(-\infty , \hat{T}]$. Moreover, 
since $q_0 ( \hat{T} ) > 0$, by lemma \ref{lem:pert} 
$\lim_{\ep \to 0} f_{\ep} (\hat{T}) = + \infty $. Hence there is 
$\ep (\omega) > 0$ such that, for $0<\ep < \ep (\omega)$,
$f_{\ep} (\hat{T}) > \cot \omega$. Now, $|q_0 (t)| = 
O (e^{\lambda_1 t} )$ as $t \to -\infty$. Hence, for all
$0<\ep < \ep (\omega)$, $\limsup_{t \to -\infty}
|f_{\ep}(t)|  \leq a(\ep)/ (|\Gamma |- a(\ep)) \leq 
a(\ep_1) / (|\Gamma| - a (\ep_1)) < \cot \omega$. Hence 
$\{ t \leq \hat{T} \ | \ f_{\ep} (t) = \cot \omega \}$
is not empty and bounded. Let $T_{\ep} $ be its smaller element
and set $r_{\ep} = |\ov{q}_{\ep } (T_{\ep}) |$. Since for all
fixed $t\leq \hat{T}$ $\lim_{\ep \to 0 } f_{\ep} (t)=0$, we must have
$\lim_{\ep \to 0} T_{\ep} = -\infty$ and, by lemma
\ref{lem:pert}, $\lim_{\ep \to 0} r_{\ep} =0$. It follows that,
provided $\ep$ is small enough, $|\ov{q}_{\ep}|$ is strictly
nondecreasing on $(-\infty , T_{\ep}]$, which yields our claim.
\end{pf}

Now, by the properties of $W_1, W_2$ and $\psi$, we have
two homoclinic solutions to (\ref{eq:pert}): $\ov{q}_{\ep}$
and $\wtilde{q}_{\ep} := - \ov{q}_{\ep}$. It remains to check,
using the previous lemma, that for $\ep$ small enough, conditions
$(H1-4)$ are satisfied. Let $\omega >0$ be small and fixed.
For $0<\ep < \ep (\omega)$, let $r_{\ep} (\omega)$ be associated to
$\omega$. 
Note that $q_{1, \ep} (-t) =
q_{1, \ep} (t)$ and $q_{2, \ep} (-t) = - q_{2,\ep}(t)$.
Hence $\ov{q}_{\ep}$ crosses the circle of center $0$ and 
radius $r_{\ep}(\omega)$ for the first time and the last time 
respectively at $\ov{\alpha} = (r_{\ep} (\omega) \cos \omega , 
r_{\ep}(\omega)  \sin \omega )$ and
$\ov{\beta} = (r_{\ep} (\omega) \cos \omega , 
-r_{\ep}(\omega)  \sin \omega )$. For 
$\wtilde{q}_{\ep}$ we have 
$\wtilde{\alpha} = (-r_{\ep} (\omega) \cos \omega , 
-r_{\ep}(\omega)  \sin \omega )$ and 
$\wtilde{\beta} = (-r_{\ep} (\omega) \cos \omega , 
r_{\ep}(\omega)  \sin \omega )$.
So it is clear that $(H1)$ is satisfied 
(with the notations of section 2, 
$\ov{\omega}_u = \omega, \ov{\omega}_s = -\omega , 
\wtilde{\omega}_u = \omega+\pi , \wtilde{\omega}_s
= - \omega +\pi$).

We have
$$
{\cal Q} := \frac{\lambda_2^2}{\lambda_1^2}
\frac{|\alpha_2| | \beta_2| + 15 (\lambda_1 / 4\lambda_2)
\Lambda r^3 }{|\alpha_1| | \beta_1|} =
\frac{\lambda_2^2}{\lambda_1^2} \Big( (\tan \omega )^2 +
\frac{15 \lambda_1}{4\lambda_2}\frac{ \Lambda r_{\ep} (\omega)}
{(\cos \omega)^2} \Big).
$$   
Now, since, when $\ep =0$, $(q_0 , 0)$ is a nondegenerate (up 
to time  translations ) homoclinic orbit, for $\ep$ small
enough $\ov{q}_{\ep}$ and $\wtilde{q}_{\ep}$ are (up to time translations)
 nondegenerate uniformly  with respect to $\ep$, and we can 
take for them a constant of nondegeneracy $C_1 (\ep)$ which is 
bounded from below by a positive constant independent of $\ep$. It follows
 that in condition $(H2)$ the second member
is bounded from below 
by some constant
independent of $\ep$. Using $\lim_{\ep \to 0} r_{\ep}(\omega)=0$,
we can derive that, provided $\omega$ is smaller than some
$\omega_0 > 0$, $(H2)$ is satisfied if $0<\ep < \ep' (\omega)$.

Since $\lim_{\ep \to 0} r_{\ep}(\omega) =0$, it is clear that
$(H4)$ is satisfied, provided $0<\ep < \ep '' (\omega)$, 
with $ \ep '' (\omega ) \in (0, \ep' (\omega)) $.

The quantity which must be greater than ${\cal Q}$ in $(H3)$ is 
$$
{\cal B} := l(\frac{\lambda_1}{\lambda_2}) 
\Big ( \frac{C_1(\ep)  {\cal M}}{36 S_2 + 28 \Lambda r_{\ep}(\omega)^2}
\Big )^{\frac{\lambda_1 }{ \lambda_2} -1} =
l(\frac{\lambda_1}{\lambda_2}) 
\Big ( \frac{C_1 (\ep)  \sin \omega }{36 \sin \omega  
+ 28 \Lambda r_{\ep}(\omega)}
\Big )^{\frac{\lambda_1 }{ \lambda_2} -1} \geq 
l(\frac{\lambda_1}{\lambda_2}) 
\Big ( \frac{C_1 (\ep)  }{36  
+ (7 \lambda_2 / 5 \lambda_1)}
\Big )^{\frac{\lambda_1 }{ \lambda_2 }-1}
$$  
for $0< \ep  < \ep''(\omega) $ and by $(H4)$.
So, if $0< \ep < \ep''(\omega)$, ${\cal B} $ is bounded from below by a 
positive constant independent of $\ep$ and $\omega$. By the 
expresion of $\cal Q$ it is clear that there is 
$\omega >0$ and $\ep_0 \in (0, \ep''(\omega))$ such that
${\cal Q } < {\cal B} $ (that is condition $(H3)$ is satisfied) 
if $0<\ep< \ep_0$. The case with $\ep < 0 $ can be dealt with in the same way.
That completes the proof of theorem \ref{th:pert}. 

\section{Relaxing the nondegeneracy assumption}\label{sec:deg}

The aim of this section is to modify the arguments of the previous 
section in order to show how to deal also with  homoclinics which are 
degenerate. For simplicity, we shall restrict
ourselves to the proof of the ``relaxed'' theorem \ref{thm:alm}, 
that is {\it thm.} \ref{thm:rela}.

\subsection{Finite dimensional reduction for 
degenerate homoclinic orbits} \label{sub:deg1}

We consider a homoclinic solution  $q$ of 
(\ref{maineq}) not necessarily nondegenerate. 
We assume that $|q(-T)|= |q(T)| =r$ and that 
$|q(t) | <r$ for all $|t|> T$.

Let $a^0 = L_{A}(c {\dot q} \chi_{[-T, T]})$, $c>0$ being
chosen such that $|a^0|_E =1$. Let $a^j$ be defined
for $1\leq j \leq p$ and satisfy $a^j = L_A e^j$, with
$e^j \in L^2 (\R)$, $e^j_{\Rs \backslash [-T,T]} \equiv 0$.
Moreover we shall assume that $(a^i , a^j ) = \delta_{ij}$. 
Let  $F= \langle a^0 , \ldots , a^p \rangle^{\bot} $ and
$\Pi$ be the orthogonal  projection on 
$F$ defined by
$$
\Pi (x) = x - \sum_{j=0}^p (x, a^j) a^j. 
$$
We shall assume that there is a constant $C'_0$ such that
$$
\forall x \in F \quad ||\Pi dS (q) x|| \geq 
C'_0 ||x||.
$$
Note that 
$Ker \ dS (q)$ is finite dimensional, so one can always
define $a^j$ enjoying the above properties.
For example, we can choose $(f^j)_{1 \leq j \leq p}$ such that
$Ker \ d S (q) \subset span \{ \dot{q} , f^1 , \ldots ,
f^p \}$ and set $a^j = L_A ( f^j \chi_{[-T, T]})$. Moreover $(a^i , a^j ) =
\delta_{ij}$ for a suitable choice of $f^1, \ldots , f^p$.

A simple  application of the Implicit function theorem leads to the following
\begin{lemma}
There is $\delta_0 >0$ and a smooth function 
$  w:  (-\delta_0 , \delta_0)^p 
\to Y$ such that :
\begin{itemize}
\item 
$ w( l^1, \ldots , l^p ) \in F$;
\item 
$S ( q + \sum_{m=1}^p l^m a^m +  w(l^1, \ldots ,l^p))
= \sum_{m=0}^p \alpha^m (l^1, \ldots, l^p) a^m $. 
\end{itemize}
Moreover $\alpha^1 (l) = \ldots = \alpha^p (l) =0 $ implies 
$\alpha^0 (l) =0$.
\end{lemma}

The last assertion is an easy consequence of the autonomy of
the system.
We shall set, for $l= (l^1,\ldots, l^p)$,
$$
q(\teta , l ) = q_{\teta} + \sum_{m=1}^p l^m a^m_{\teta} + w(l)_{\teta}.
$$
\begin{remark}
By the definition of $q(\teta , l)$ and the autonomy of the system, 
we  have  $q(\teta,l) = q(0,l)_{\teta}$ and 
$S(q(\teta ,l))= \sum_{m=0}^p \alpha^m (l) a^m_{\teta}.$
Moreover we have  the estimate
\be \label{eq:q(teta)}
||q(\teta,l)- q_{\teta}|| \leq \ov{C} \max_m|l^m| 
\leq \ov{C} \delta_0 ,
\ee
where $\ov{C}$ depends on $C'_0$.
We can derive that $q(\teta ,l)$, by the equation it satisfies,
belongs to $X$ (provided $\delta_0$ has been chosen small enough,
namely $\ov{C} \delta_0 < \min \{\rho_0, 2/\Lambda \})$.
\end{remark}
Set $ {\cal G} (l) = 
f ( q (0, l))= f(q(\teta ,l))$.
By the properties of $w(l)$ we can easily prove the following lemma
\begin{lemma} \label{lem:redg}
$(\partial {\cal G} / \partial l^m) (l) 
= \alpha^m (l)$. As a consequence
${\cal G}'(l)=0$ iff $q(\teta , l)$ is (for all $\teta$) a homoclinic
solution of (\ref{maineq}).
\end{lemma}
\begin{remark} \label{rem:deg2}
Reciprocally, there is a neighborhood $U$ of $q$ in 
$Y$ such that all the homoclinic solutions in $U$ 
correspond to critical points of ${\cal G}$.
\end{remark}

As in section \ref{sec:finite} 
the following estimate holds (provided $R$ is large 
enough):  there is $C'_1 > 0$ such that for all $x \in F$
\be \label{eq:S'}
\max \Big( ||dS (q) x -  \sum_{m=0}^p \eta^m a^m || ,
R |(a^0,x)| , \ldots , R|(a^p,x)| \Big) \geq C'_1 
\max ( ||x|| , |\eta^0|, \ldots , |\eta^p| ).
\ee
Now assume that we have two distinct homoclinic solutions 
$\ov{q}$ and $\wtilde{q}$, and finite dimensional
spaces $\wtilde{\ov{F}} = \langle 
\wtilde{\ov{a}}^0, \ldots , \wtilde{\ov{a}}^p \rangle^{\bot}$,
with the above properties.

Then, provided $\delta_0 $ is small enough, 
for $j=(j_1, \ldots , j_k) \in \{ 0,1 \}^k$,
for $ \Theta = (\teta_1 , \ldots , \teta_k) \in \R^k$, 
$L= (l_1, \ldots , l_k) \in ((-\delta_0 ,\delta_0)^p)^k$, with 
$\ov{d} > 2/ \lambda_2 $,  we can build $Q(\Theta , L)$ in 
the same way as we  built 
$ Q_{\Theta}$ in subsection \ref{sec:nat}, just substituting
 $q(\teta_i , l_i)$ to $q_{\teta_i}$ in   the construction.
Note that we keep the same $T_i$ satisfying $|\wtilde{\ov{q}} (T_i)|=
r$ so that in this case our boundary value problems may connect two
points with different norms (however we know that these norms
are $\leq r + \ov{C} \delta_0$).

In the spirit of lemma \ref{lem:grad}, one can get
\be \label{eq:useful}
||S(Q(\Theta , L))- \sum_{i=1}^k \sum_{m=1}^p
\alpha^m (l_i) a_i^m || \leq K (r + |L| ) e^{-\lambda_2 \ov{d}}.
\ee
Here is the equivalent of lemma \ref{lem:const}. 
It is not so specific, 
but it is enough to prove the equivalent of theorem \ref{thm:alm}.
We shall use the notations $a^m_i = \ov{a}^m_{\teta_i}$ 
if $j_i =0$, 
$a^m_i = \wtilde{a}^m_{\teta_i}$  if $j_i =1$.
\begin{lemma} \label{lem:deg3}
There are  ${\ov D}_1$ and $\delta_1 >0$ which depend 
on $C'_1, \Lambda, \ov{\Lambda}$ and
 there exist  $\ov{w}$,  function of
$\Theta = (\teta_1 , \ldots , \teta_k)$ and of 
$L= (l_1 ,\ldots , l_k)$, defined for 
$\ov{d} \geq {\ov D}_1$ and $ l_i \in ( - \delta_1 , \delta_1)^p$,
such that $(\Theta , L) \mapsto Q(\Theta , L) +
\ov{w} (\Theta , L)$ is smooth and
\begin{itemize}
\item $\ov{w} (\Theta , L) \in \cap_{i=1}^k
\langle a^0_{\teta_i}, \ldots , a^p_{\teta_i} \rangle^{\bot} $;
\item 
$S ( Q(\Theta , L) + \ov{w}(\Theta , L) ) = \sum_{i=1}^k 
\sum_{m=0}^p \alpha_i^m (\Theta , L) a^m_{\teta_i} $.
\end{itemize}
Moreover
$$
||\ov{w}|| \leq K_1 (r + |L|)  e^{-\lambda_2 \ov{d}},
$$
where $|L| = \max_{i,m} |l_i^m|$ and 
$K_1$ depends only on $C'_1$ and $\ov{\Lambda}$,$\Lambda$.
\end{lemma}
The proof can be carried out in the same way as in the 
nondegenerate case.
Set 
$$
g (\Theta , L) = f ( Q( \Theta, L) + \ov{w} (\Theta , L)) .
$$
We have, by lemmas \ref{lem:one} and \ref{lem:deg3}
\begin{lemma}
Every critical point of $g $ gives rise to a $k$-bump 
homoclinic solution to the system, provided $\delta_1$ and
$r$ have been chosen small enough.
\end{lemma}
Here again the proof does not differ from the one given 
in the nondegenerate case. 

Finally, using the notations $\ov{{\cal G}}(l)=
f( \ov{q} (0,l)) $ and $\wtilde{{\cal G}}(l) = 
f ( \wtilde{q} (0,l))$,  we get:
\begin{lemma} \label{lem:deg5}
For all
 $\Theta = ( \th_1, \ldots, \th_k )$ 
with $\ov{d} > {\ov D}_1 $ and $L$ with
$|L| \leq \delta_1 $ 
\begin{eqnarray*} 
 &  (i) & \frac{\partial}{\partial d_i} g ( \Theta, L) =
\frac{\partial}{\partial d_i} e
\Big( (Q(\Theta, L) + \ov{w} (\Theta, L) )( s_i), 
(Q(\Theta, L) + \ov{w} (\Theta, L) )( u_{i+1}),
d_i \Big) \\
&(ii)  & \Big| \frac{\partial g }{\partial l_i^m } - 
\partial_{l^m}  {\cal G}_i (l_i^1, \ldots , l_i^p) \Big| 
\leq K_2 (r + |L|) e^{-\lambda_2 \ov{d}},
\end{eqnarray*}
where $K_2$ depends only on $K_1$ and $\Lambda$, $\ov{\Lambda}$ 
and ${\cal G}_i = \ov{{\cal G}} $ if $j_i=0$, 
${\cal G}_i = \wtilde{\cal G} $ if $j_i =1$.
\end{lemma}
\begin{pf}
$(i)$ is proved exactly in the same way as lemma 
\ref{lem:funct}.
For $(ii)$, write
$$
\frac{\partial}{\partial l_i^m} g = 
\Big( S(Q(\Theta ,L)+ \ov{w} (\Theta, L) ), \frac{\partial
Q (\Theta,L)}{\partial l_i^m} + \frac{\partial \ov{w}
(\Theta,L) }{\partial l_i^m} \Big) . 
$$
We have by lemma \ref{lem:deg3}
$S(Q(\Theta ,L)+ \ov{w} (\Theta, L) ) = \sum_{i,m}
\alpha_i^m (\Theta,L) a_i^m . $
Since $\ov{w} (\Theta,L) \in 
\cap_{1\leq i\leq k} \cap_{0\leq m \leq p} 
\langle a_i^m \rangle^{\bot}$,
$\partial \ov{w}/\partial l_i^m (\Theta ,L) $ belongs to the same
space. Hence, since supp $ \partial Q (\Theta , L)/
\partial l_i^m \subset [\teta_{i-1}+ T_{i-1}, \teta_{i+1}
- T_{i+1} ]$  and supp $(-\ddot{a}_i^m + A a_i^m) 
\subset [\teta_i -T_i, \teta_i +T_i ]$,
$ (\partial  g)( \partial l_i^m ) = $
$$ 
\Big( \sum_{n,q}
\alpha_n^q (\Theta,L) a_n^q , \frac{\partial}{\partial l_i^m}
Q (\Theta , L) \Big) 
= \sum_{q=0}^p \alpha_i^q (\Theta , L) \Big( a_i^q , \frac{\partial}
{\partial l_i^m} q^i ( \teta_i , l_i) \Big)
= \sum_{q=0}^p \alpha_i^q (\Theta , L)
(a_i^q , a_i^m )= \alpha_i^m (\Theta , L) .
$$
We have used there the definition of $q(\teta , l)$,
$\partial q (\teta_i , l_i )/ \partial l_i^m =
a^m_i + (\partial_{l^m} w_i(l_i))_{\teta_i}$ and 
$(\partial_{l^i} w_i(l_i))_{\teta_i} \in
\langle a_i^0, 
\ldots , a_i^p  \rangle^{\bot}$, where
$w_i = \ov{w}$ (resp. $w_i = \wtilde{w}$) if 
$j_i =0$ ( resp. $j_i =1$).
We get 
\begin{eqnarray*}
\frac{\partial }{\partial l_i^m} g  & = & 
(S(Q(\Theta ,L))+ \ov{w}(\Theta , L)), a_i^m )\\
&=& (S(Q(\Theta ,L),a_i^m ) + O( (r+|L|) e^{-\lambda_2 \ov{d}})
\\
&=& \alpha_i^m (l_i) + O( (r+|L|) e^{-\lambda_2 \ov{d}}) \\
&=& \partial_{l^m} { \cal G }_i (l_i) + O( (r+|L|) e^{-\lambda_2 \ov{d}}),
\end{eqnarray*}
where we have used lemma \ref{lem:deg3} and 
(\ref{eq:useful}) in the second and the third line 
respectively. This is exactly (ii).
\end{pf}

Now we state a corollary of lemma 
\ref{lem:deg5} (i)  which is  got from  a simplified
version of lemma \ref{lem:ener}.

\begin{corollary} \label{cor:deg1}
For all $(\Theta, L)$ with $\ov{d} > {\ov D}_1 $ and
$|L| \leq \delta_1 $ 
\be \label{eq:estfo}
\Big| \frac{\partial}{\partial d_i} g (\Theta , L) -
2 \Big( \lambda_1^2 \alpha^{i+1}_1 \beta^i_1 e^{-\lambda_1 d_i} +
\lambda_2^2  \alpha^{i+1}_2 \beta^i_2 e^{-\lambda_2 d_i} \Big) \Big|
\leq K_4 e^{-\lambda_2 d_i} \Big( (r+|L|)^2 e^{-\lambda_2 \ov{d}} 
+ (r+|L|)^3 + r |L| + |L|^2 \Big),
\ee
where $K_4$ depends only on $C'_1$ and $\Lambda$, $\ov{\Lambda}$.
\end{corollary}
\begin{pf}
We omit the details of the proof. It is a simple 
consequence of lemmas \ref{lem:deg3},\ref{lem:deg5}
 and \ref{lem:ener}.
\end{pf}

\subsection{Topological nondegeneracy} \label{sub:top}

Let $q$ be a (possibly degenerate) {\it isolated} homoclinic
solution of (\ref{maineq}). ``Isolated'' means here that there is 
a neighborhood $U$ of $q$ in $Y$ such that all the 
homoclinics which belong to $U$ are translates of $q$.
Let $a= L_A (c\dot{q} \chi_{[-T, T]}) \in X $, where $c>0$ is chosen so 
that $|a|_E =1$. Let
$\hat{F} = a^{\bot}$ and  
 $\hat{\Pi} : Y
\to \hat{F} $ be the projection  defined by
$$
\hat{\Pi} (x)= x - (x,a) a. 
$$
Consider $G_* : \hat{F} \to \hat{F} $, defined by
$$
G_* (x)=\hat{\Pi} S(q+x) =
\hat{\Pi} \Big[ (q+x)- L_A \Big( \nabla W (q+x) - \psi (q+x){\cal J}
 (\dot{q}+\dot{x}) \Big) \Big] 
 = \hat{\Pi} ( S (x) - K_q (x) ),
$$
where
$$
K_q (x) = L_A \Big[ \Big( \nabla W (q+x) -
\nabla W (q) - \nabla W (x) \Big) -
\Big( (\psi (q+x)- \psi(q) ) {\cal J}  \dot{q} +
(\psi (q+x) - \psi(x)) {\cal J} \dot{x} \Big) \Big].  
$$ 
We have $K_q (0)=0$. In addition,
it is easy to see that $K_q$ sends $Y$ into
$E$ and that it is compact.

Note also that  there is $\rho > 0$ 
such that $ \hat{\Pi} \circ S : \hat{F} \to 
\hat{F}$ is a diffeomorphism 
from ${\hat B} (\rho)$ onto a neighborhood of $0$ in
$\hat{F}$ containing ${\hat B}(\rho / 2)$. Let $\delta > 0$ satisfy
$\hat{\Pi} K_q ({\hat B}(\delta)) \subset {\hat B}(\delta/2) $ and
let $\hat{G} : {\hat B}(\delta) \to \hat{F}$ be defined by
\be \label{eq:tildeG}
\hat{G}(x) = x - (\hat{\Pi} \circ S )^{-1} 
\hat{\Pi} K_q (x) :=  x - \hat{K}_q (x).
\ee
We have $\hat{K}_q (0) =0$, and $\hat{K}_q (Y)
\subset E$. Hence all the zeros of $\hat{G}$ must belong
to $E$ and thus be homoclinic solutions to
the system. Now $q$ being an isolated homoclinic, $0$ is 
an isolated zero of $\hat{G}$.

Moreover
$\hat{K}_q$ is a compact operator. We can now introduce
the following definition : 
\begin{definition}\label{eq:topdeg}
We shall say that $q$ is a ``topologically nondegenerate'' 
homoclinic if there is $0< \nu  \leq \delta$
such that $\deg (\hat{G}, \hat{B} (\nu) , 0) \neq 0$ and $\hat{G}$ has no
zero in $\ov{{\hat B}(\nu)}\backslash \{ 0 \}$.
\end{definition}

\begin{remark}
We could prove without difficulty that this definition 
is independent of the choice of $a$ satisfying 
$(a, \dot{q}) \neq 0$.
\end{remark}

As a consequence  to the hyperbolicity of the equilibrium, the 
solutions of the linearised system about $q$ which belong
to $Y$ must belong to $X$. Therefore, if $q$ is a 
nondegenerate homoclinic, $dG_*(0)$   and hence
$d\hat{G}(0)$ are injective. By 
(\ref{eq:tildeG})  $\hat{G}$ is then a 
local diffeomorphism about $0$ and the 
property of ``topological nondegeneracy ''
defined above is satisfied.

\begin{remark}
We point out that in certain cases, one can say 
that a variationally obtained isolated  homoclinic is topologically
nondegenerate. For instance, an isolated local minimum for 
$f$, or, under further conditions, an  isolated mountain-pass critical point
correspond to topologically nondegenerate homoclinics
(see \cite{A}, \cite{HH} and \cite{BS}).
\end{remark}

Now consider as in subsection \ref{sub:deg1} 
$a^0 ,\ldots , a^p$ which satisfy the properties given in this 
subsection. We can then define the function ${\cal G}$ on some 
$(-\delta_0 , \delta_0 )^p$. We shall prove

\begin{lemma}
Assume that the homoclinic $q$ is isolated and topologically 
nondegenerate. Then
$0$ is an isolated critical point of ${\cal G}$. Moreover
there is $\mu \in (0, \delta_0) $ such that 
$\deg ({\cal G}', (-\mu , \mu)^p , 0) \neq 0$ and ${\cal G}'$ has no zero in
$(-\mu , \mu)^p$, except $0$.
\end{lemma}
\begin{pf}
By remark \ref{rem:deg2}
since  $q$ is an isolated homoclinic, $0$ is an 
isolated critical point of ${\cal G}$  and there is $\mu_1$ such 
that ${\cal G}'$ has no zero in $(-\mu_1 , \mu_1)^p$ except $0$.

Let $\hat{G} $ be the function defined above,
associated to $a_0$.
By the topological nondegeneracy property of $q$, there is some
$\nu > 0$ such that $\deg(\hat{G}, \hat{B} (\rho),0)\neq 0$ 
for all $0<\rho< \nu$. 
Let $B_{\nu_2} = \{ x \in F \ | \ || x || \leq \nu_2  \}.$
Consider for $\delta_2, \nu_2$ small enough 
$ \varphi : (-\delta_2 , \delta_2)^p \times B_{\nu_2} \to 
{\hat B}(\nu) $ 
which assigns to
$ (l,y) \to $ $\varphi(l,y)= 
q(0 ,l ) + y - q_0 =  \sum_{j=1}^p l^j a^j + w(l) + y $.
This is clearly a diffeomorphism from 
$ (-\delta_2 , \delta_2)^p 
\times B_{ \nu_2} $ onto some neighborhood of $0$
in $\hat{F}$ included in $\hat{B} (\nu)$.
Now, let $\xi : (-\delta_2 , \delta_2)^p 
\times  B_{\nu_2} \to \hat{F}$ be defined by
$$
\xi (l, y) = \hat{G} (\varphi (l, y))=
\varphi (l,y) + \hat{K}_q ( \varphi (l,y)).
$$
Since the degree is invariant under a diffeomorphism there results that 
$$
\deg (\xi, (-\delta_2 , \delta_2)^p 
\times  B_{\nu_2} , 0) =
\deg(\hat{G}, \hat{B} (\nu), 0) \neq 0 . 
$$
We decompose $ G (\varphi (l,y)) = A(l,y) + \sum_{j=1}^p u^j(l,y)a^j$,
where $ (A(l,y) , a^j) =0$.
For $t\in [0,1]$, define $\xi_t$   by
$$
\xi_t (l,y) = (\Pi \circ S )^{-1} \Big(  A((1-t)l, y) 
+ \sum_{j=1}^p u^j (l, (1-t)y)a^j \Big).
$$ 
$\xi_t$ has the form $I-K_t$, where $K_t$ is compact.
Moreover, for all $l \in (-\delta_2,
\delta_2 )^p$, $ A(l, y)=0 $ iff $y=0$. Hence 
$\xi_t (l,y) =0$ iff $y=0$ and $ u^j (l,0)=0$. Now
$u^j (l,0) = \alpha^j (l) $, and since $q$ is an 
isolated homoclinic, it vanishes at no other point in 
$(-\delta_2 , \delta_2)^p $ than $0$. Therefore
$$
\deg (\xi , (-\delta_2 , \delta_2)^p 
\times  B_{\nu_2},0) =
\deg (\xi_1 , (-\delta_2 , \delta_2)^p 
\times  B_{\nu_2} , 0 ).
$$
Now, $A(0,\cdot ) $ is a  diffeomorphism from 
$B_{\nu_2} $ to a neighborhood of $0$
in $ F $. Let $\Psi$ be 
defined on $\hat{B} (\nu_2)$ by 
$$
\Psi (y + \sum_{j=1}^p r_j a^j) = A(0,y) + \sum_{j=1}^p r_j a^j.
$$
Since $\Psi$ is a diffeomorphism from $B(\nu_2)$ to a neighborhood
of $0$ in $\hat{F}$ and we have 
$$
\xi_1 (l,y) = (\Pi\circ S)^{-1} (
\Psi(y+ \sum_{j=1}^p \alpha^j (l) a^j)),
$$
with $(\Pi \circ S)^{-1} (\Psi(0))=0$. Hence, setting
$\tau (l,y) = \sum_{j=1}^p \alpha^j (l) a^j + y$ we get 
$$
|\deg ( \xi_1 , (-\delta_2 , \delta_2)^p 
\times  B_{\nu_2}, 0)| = 
|\deg (\tau , (-\delta_2 , \delta_2)^p 
\times  B_{\nu_2} , 0)| = | \deg ((\alpha^1,
\ldots , \alpha^p ) ,
(-\delta_2 , \delta_2 )^p, 0)|. 
$$ 
Using lemma \ref{lem:redg}, we get the result
with $\mu = \min \{ \delta_2, \mu_1 \}$.
\end{pf}
\subsection{Relaxed theorem for system $(S_{\ep})$} \label{sub:relth}

We shall prove

\begin{theorem}\label{thm:rela}
Assume that the hypotheses of theorem \ref{thm:alm} hold, with the 
nondegeneracy condition replaced by topological nondegeneracy.
Then the same conclusion holds.
\end{theorem}

\begin{pf}
$\ov{q}$ and $\wtilde{q}$ may be degenerate. Since the configuration
space is $\R^2$, $W^u (0) $ and $W^s (0)$ are 2-dimensional.
Therefore  dim $Ker \ dS (\wtilde{\ov{q}}) $ cannot exceed $2$.
Hence in the construction of subsection \ref{sub:deg1} the spaces 
associated to $\ov{q}$ and $\wtilde{q}$, 
span $ \{ \wtilde{\ov{a}}^0, \wtilde{\ov{a}}^1 \}$,  
are two dimensional.  

The topological nondegeneracy of $\ov{q}$ and
$\wtilde{q}$ implies that, for $\ep$ small enough, 
$(S_{\ep})$ has two homoclinics $\ov{q}_{\ep} $ and
$\wtilde{q}_{\ep}$such that 
\be \label{eq:dqep}
\lim_{\ep \to 0} \max \{ |\ov{q}_{\ep} - \ov{q}|_{\infty},
|\wtilde{q}_{\ep} - \wtilde{q}|_{\infty} \} =0.
\ee
 Moreover, there is 
a constant $C'_2$ independent of $\ep$ small such 
that $|| dS_{\ep} ( \wtilde{ \ov{q}_{\ep} } ) x || \geq 
C'_2 ||x||$ for $(x, \wtilde{\ov a}^1) = (x, \wtilde{\ov a}^2)=0 $.
As a consequence we may define $ \wtilde{\ov{q}}_{\ep} 
(\teta ,l)$ for $|l| \leq \mu_0$. Writing
$\wtilde {\ov{ {\cal G} }}_{\ep}(l) = 
f_{\ep} (\wtilde{ \ov{q} }_{\ep} (\teta , l)) $ ,
we know that $|{\cal G}'_{\ep}- {\cal G}'|_{{L^{\infty}(-\mu_0 ,\mu_0)}}
\to 0$ as $\ep \to 0$. Note that it may occur that
${\cal G}_{\ep}$ has a sequence of critical points 
converging to $0$. So we cannot say that $\ov{q}_\ep,
\wtilde{q}_{\ep}$ are 
 isolated homoclinic. However we know that, for 
all $|\mu | \leq \mu_0$ there is $\ov{\ep}$ such that,
for $ |\ep| < \ov{\ep}$, all the critical points of ${\cal G}_{\ep}$
belong to $(-\mu , \mu)$ and 
$\deg({\cal G}'_{\ep}, (-\mu ,\mu),0 )\neq 0$.

Given $k$ and $j=(j_1, \ldots, j_k) \in \{ 0,1 \}^k$, 
we can construct $Q_{\ep}(\Theta ,L)$ as well
as $\ov{w}_{\ep} (\Teta , L)$ for all 
$L=(l_1, \ldots ,l_k) \in (-\mu_0 , \mu_0)^k $ and 
for all $ \Theta = (\teta_1, \ldots, \teta_k) \in \R^k$ with
$ \ov{d} \geq \ov{D}_1$.
We can as well define $g_{\ep} ( \Theta , L) $.  

By the properties of the system near the equilibrium,
 there are $r_4>0$ and $\ep_2 >0$ such that,  for $0<r < r_4$
and $|\ep| < \ep_2$,
the trajectory of $\ov{q}_{\ep}$ (resp. 
$\wtilde{q}_{\ep}$) crosses the circle of radius
$r$ at two points only:  
$\ov{\alpha}_{\ep} (r)$ and $\ov{\beta}_{\ep}(r)$ 
(resp. $\wtilde{\alpha}_{\ep} (r)$ and $\wtilde{\beta}_{\ep}(r)$).
Moreover, 
$$
\lim_{r \to 0} \frac{\ov{\alpha}_0 (r)}{r }=
(\cos \ov{\om}_u , \sin \ov{\om}_u ) ; \quad
\lim_{r \to 0} \frac{\ov{\beta}_0 (r)}{r }=
(\cos \ov{\om}_s , \sin \ov{\om}_s )
$$
and for all $r \in (0,r_4)$, from (\ref{eq:dqep}),
\be \label{eq:eglim1}
\lim_{\ep \to 0} \ov{\alpha}_{\ep}(r)
=\ov{\alpha}_{0}(r) ; \quad
\lim_{\ep \to 0} \ov{\beta}_{\ep} (r)
=\ov{\beta}_{0} (r).
\ee
We have similar properties for $\wtilde{q}, \wtilde{q}_{\ep}$.
Let us define
$ P^i_1 = \cos \om^i_u   \cos \om^{i+1}_s,$
$P^i_2 = \sin  \om^i_u   \sin \om^{i+1}_s$.
Lemma \ref{lem:deg5} and Corollary 
\ref{cor:deg1} hold ( with $g$ replaced by $g_{\ep}$),
and $\alpha, \beta$ replaced by $\alpha_{\ep}, \beta_{\ep}$
we have 
\be \label{eq:thbis}
\Big| \frac{\partial}{\partial d_i} g_{\ep} -
2 r^2 \Big( \lambda_1^2  P^i_{1,\ep} e^{-\lambda_1 d_i} + 
\lambda_2^2 P^i_{2,\ep} e^{-\lambda_2 d_i} \Big) \Big| \leq K_5 
(r^3 + r \nu +\nu^2 + r^2 e^{-\lambda_2 \ov{d}})e^{-\lambda_2 d_i}, 
\ee 
where $K_5$ is independent of $\ep $ small and 
$ |P^i_{1,\ep} - P^i_{1}| $, 
$ |P^i_{2,\ep} - P^i_2| \leq $ $ r(\ep)$
$ r(\ep) \to 0$   as $ \ep \to 0.$
We recall that $\lambda_1 = \lambda + |\ep |$ and
$\lambda_2 = \lambda - |\ep |$.
We know that $P^i_1 P^i_2 < 0$ and that there exists $\delta >0$ such that, 
for all $i$,
$$
\frac{|P^i_1|}{|P^i_2|} \geq 1+ \delta > 1 \quad ,
\quad |P^i_1|,|P^i_2| \geq \delta > 0. 
$$

First choose $0<r < r_5$ and $\nu_1 > 0$ such that 
$$
\frac{|P^i_1|}{|P^i_2|+ K_5 (r + \nu_1 /r + \nu_1^2 /r^2)} 
\geq 1+ \delta /2  \quad ,
\quad |P^i_{1,2}-K_5 (r + \nu_1 /r + \nu_1^2 /r^2) | \geq \delta /2 .
$$
Then  set
$$
D = \frac{1}{2\ep} \min_i \ln \left(   
\frac{|P^i_1|}{|P^i_2|+ K_5 (r + \nu_1 /r + \nu_1^2 /r^2)} \right) - 1 \ ; \
J  = \frac{1}{2\ep} \max_i \ln \left(   
\frac{|P^i_1|}{|P^i_2|- K_5 (r + \nu_1 /r + \nu_1^2 /r^2)} \right) + 1. 
$$
By (\ref{eq:thbis}) it is easy to see that 
$\lim_{\ep \to 0} D = + \infty$ and that  there is $\ep_3>0$ such that 
for all $|\ep| < \ep_3$, for all $(d_1 , \ldots , 
d_k) \in (D,J)^k$, for all $L = (l_1 , \ldots , l_k) \in (-\nu_1 , \nu_1)^k$,
$$
d_i =D \Rightarrow {\rm sign} \Big( \frac{\partial}{\partial d_i} 
g_{\ep} \Big) = - {\rm sign} \ P^i_2  \quad  , \quad 
 d_i =J \Rightarrow {\rm sign} \Big( \frac{\partial}{\partial d_i} 
g_{\ep} \Big) =  {\rm sign} \ P^i_2 
$$ 
and 
$$
|L| = \nu_1 \rightarrow 
|\partial_{l^m} g_{\ep} - \partial_{l^m} {\cal G}_m | < 
|\partial_{l^m} {\cal G}_m|/2. 
$$
Now arguing as in the proof of theorem 
\ref{thm:alm} and using that
$\deg ( \partial_{l^m} {\cal G}_m , (-\nu_1 , \nu_1 ) ,0) \neq 0 $,
we get 
$$
\deg (d g_{\ep}, (D,J)^k \times (-\nu_1 , \nu_1 )^k ,0 )
\neq 0,
$$ 
which implies the desired result.
\end{pf}

\section{Dynamical consequences}\label{sec:dyn}

A family of multibump solutions like the ones of theorem \ref{thm:main2} 
ensures the positivity of the topological entropy
at  the zero energy level ${\cal E}^{-1}(0)$, see also \cite{BS}
and \cite{S2}. 
We denote by $\Phi ( t , x) \in \R^4$ with $x = ( q, \dot{q})$ 
the flow associated to (\ref{maineq}).
The definition of the  topological entropy is the following:
$$
h_{top} = \sup_{R, e > 0} 
\left( \limsup_{t \to + \infty } \frac{\log s(t,e,R) }{t} \right)    
$$
where
\begin{eqnarray*}
s(t,e,R) &=& \max \{ {\rm Card}(\wtilde{E}) \ | \  \forall \tau \in [0,t] 
\ : \ 
\Phi (\tau , \wtilde{E} ) \subset B_4 (0,R), \\
& & \forall x \neq y \in \wtilde{E}, 
\exists \tau \in [0,t] \ : \ |\Phi(\tau, x) -\Phi(\tau, y)| \geq e > 0 \}.  
\end{eqnarray*}

We formulate the following corollary of theorem \ref{thm:main2}.
\begin{theorem}\label{thm:infin}
Assume $(W1)$,$(P1)$,$(v1)$,$(S1-2)$ and $(H1-4)$. 
There exist $0<D<J$ such that for every sequence 
$j \in \{ 0,1 \}^{\Z}$ there is $\Th \in {\R}^{\Z} $ 
 with $d_i \in ( D, J )$ 
and a solution $x_j$ of system (\ref{maineq}) such that
\begin{itemize}
\item if $j_i =0$ then on the interval $[ \th_i - {\ov T} , \th_i + \ov{T}]$ 
$$
|x_j (t) - {\ov q}( t - \th_i ) | \leq 
\frac{r}{8}
\min (|\cos \ov{\om}_{u,s}|, |\cos \wtilde{\om}_{u,s}|,
|\sin \ov{\om}_{u,s}|, |\sin \wtilde{\om}_{u,s}|), 
$$
\item if $j_i =1$ then on the interval $[ \th_i - {\wtilde T} , \th_i + 
\wtilde{T}]$ 
$$
|x_j (t) - {\wtilde q}( t - \th_i ) | \leq 
\frac{r}{8}
\min (|\cos \ov{\om}_{u,s}|, |\cos \wtilde{\om}_{u,s}|,
|\sin \ov{\om}_{u,s}|, |\sin \wtilde{\om}_{u,s}|), 
 $$
\item
Outside $(\cup_{j_i =0} [\theta_i - \ov{T} , 
\theta_i + \ov{T}]) \cup (\cup_{j_i =1} [\theta_i - \wtilde{T},
\theta_i + \wtilde{T}] ) $, $|x_i (t) | \leq 2r$.
\end{itemize}
\end{theorem}
 
We want to estimate $s(t^*_k,e,R^*)$ with 
$R^* = \max \{ |(\ov{q}, \dot{\ov{q}})|_{\infty}, 
|(\wtilde{q}, \dot{\wtilde{q}})|_{\infty} \}+r$,
$t^* = (k-1)( 2 \max \{ \ov{T}, \wtilde{T} \} + J)$.
Let  
$
A = \{ j \in \{ 0,1 \}^{\Z} \ | \ 
j_i = 0 \ {\rm for} \ i < 0 \ {\rm and } \ 
j_i = 0 \ {\rm for}  \ i \geq k \}
$
Associate to $j \in A$  
a solution $x_j$ given by theorem \ref{thm:infin}, for which we
may assume, by the autonomy of the system, that $|x_j (0)|=r$
and $|x_j (t)| < r, \ \forall t<0$. Consider
$\wtilde{E} = \{ (x_j, {\dot x}_j) (0)  \ | \  j \in A \}$.

As an easy consequence of  theorem
\ref{thm:infin} and Hypotheses ($H1-4$),  
there is $e>0$ such that, if $j \neq j' \in A$,
then there is $\tau \in [0, t^*_k]$ such that
$|\Phi(\tau, x_j (0) , \dot{x}_j (0) )-
\Phi(\tau, x_{j'} (0) , \dot{x}_{j'} (0) ) | \geq e$.

Hence $s(t_k^*,e ,R^*) \geq {\rm Card} \wtilde E = 2^k$.
Since $ 1 / t^* > 1 / (k-1) ( 2 \max \{ \wtilde{T}, \ov{T} \} + J )$ 
we finally deduce that $$ h_{top}^0 > 
\frac{\log 2}{2 \max \{ \wtilde{T}, \ov{T} \} + J}.$$

\begin{remark}\label{rem:dimn}
The above results could  be generalized to systems 
like (\ref{maineq}) with $q \in \R^n$, where $0$ is a
hyperbolic equilibrium of characteristic exponents $\pm \lambda_i$,
with $\lambda_1 > \lambda_2 \geq \ldots \geq \lambda_n$.
The finite dimensional reduction could be performed  
in the same way and  we would have  to impose conditions  
similar to $(H1-4)$ in order to be able to get multibump
homoclinic solutions.
\end{remark}

\section{Appendix}

We shall assume everywhere that $d \geq 2/\lambda_2$
and that $0 < r < \rho_0 /2 $.

\begin{pfn}{\sc of lemma} \ref{lem:homo}.
We perform the proof for $q^{+}_h$.

The proof of the uniqueness assertion is left to the reader (in
fact there is uniqueness also in the class of the functions 
$[0, +\infty) \to B_{2r}$.

For the existence proof let us define the Banach space:
$$
Z_1= \Big\{ g \in W^{1, \infty}  [ 0 , + \infty )  \ | \  
|g(t)| e^{\lambda_2 t},
| {\dot g}(t)|  e^{\lambda_2 t} \in L^{\infty} \Big\} ,
$$
endowed with the  norm 
$$ 
||g||_1 =  
\max \Big\{ \sup_{t \in [0, +\infty ) } |g(t)| e^{\lambda_2 t}, 
\sup_{t \in [0, +\infty )}
\frac{1}{\lambda_2} |{\dot g}(t)| e^{\lambda_2 t} \Big\} .
$$
We call $q_l = q_{h,L}^+  = e^{-t \sqrt{A} } \beta$. 
Our problem is equivalent to finding a fixed point in $Z_1$ of 
\be\label{eq:fix} 
 {\cal F}(x) := {\cal L} (\nabla W(q_l
 +x) - \psi(q_l +x) {\cal J} ({\dot q}_l +
\dot{x})),
\ee
where
 ${\cal L}$ is the  operator which 
assigns to $h$ the unique solution $ u ={\cal L} h $ of the problem:
$$
 - {\ddot u } + A u = h
\quad {\rm with } \quad u(0)= 0 \ {\rm and } \ 
\lim_{t \to +\infty}  u(t) = 0.
$$
An explicit computation shows that:
\be \label{eq:s}
{\cal L} h (t) = 
\frac{1}{2} ( \sqrt{A})^{-1}\left[ 
\int_{t}^{+\infty} 
\Big( e^{(t-s) \sqrt{A} }-e^{-(t+s) \sqrt{A} } \Big) h(s) ds +  
\int_{0}^t 
\Big( e^{(s-t) \sqrt{A} } -e^{-(s+t) \sqrt{A} } \Big) h(s) ds \right]
\ee 
and
$$
\frac{d}{dt} {\cal L} h (t)=
\frac{1}{2} 
\int_{t}^{+\infty}
\Big( e^{(t-s) \sqrt{A} }+e^{-(t+s) \sqrt{A} } \Big) h(s) ds - 
\frac{1}{2}  
\int_{0}^t 
\Big( e^{(s-t) \sqrt{A} } -e^{-(s+t) \sqrt{A} } \Big) h(s) ds.
$$
Call $B_{\delta}= \{ x \in Z_1 \ | \ ||x||_1 < \delta \}. $
We want to solve (\ref{eq:fix}) by means of the 
contraction mapping theorem in $B_{\delta}$.
So we want to find, for $r$ small enough, $\delta <r$ small enough
in such a way that: 
\\[1mm]
{\bf (i)} $\ov{{\cal F} ( B_{\delta} )} \subset B_{\delta}\quad $;
{\bf (ii)} ${\cal F}$ is a contraction on $\ov{ B_{\delta}} $.
\\[1mm]
\indent
Assume that $x \in B_{\delta}$. Then, by $(W1)$ and $(P1)$ 
$$
\Big| \nabla W (q_l +x) -  
\psi(q_l +x){\cal J} (\dot{q}_l +\dot{x} ) \Big| \leq 
\Big( \frac{1}{2} L_1 (r + \delta)^2 + 
L_2 (r+\delta) (\lambda_1 r + \lambda_2 \delta) \Big) e^{-2 \lambda_2 t}
\leq \frac{\Lambda}{2}\lambda_2^2  (r+\delta)^2  e^{-2 \lambda_2 t}.
$$
Hence, for $x \in B_{\delta}$,
$$
|{\cal F}(x)(t)|, \frac{1}{\lambda_2} \Big| \frac{d}{d t} {\cal F}(x) (t) \Big|
 \leq \frac{\Lambda}{4} \lambda_2  (r + \delta)^2   
\Big[ \int_t^{+\infty} e^{\lambda_2 (t-3s)} \ ds +
\int_{t}^{+\infty} e^{- \lambda_2 (t+3s)} \; ds 
+ \int_0^t e^{- \lambda_2 (t+s)} \; ds \Big].
$$
We get 
\be \label{eq:prim}
||{\cal F}(x)||_1 \leq  \frac{\Lambda}{4} (r+\delta)^2.
\ee
On the other side, elementary estimates give that
by $(W1)$ and $(P1)$
$$
\Big| \nabla W (q_l +x') -  \psi(q_l +x'){\cal J}(\dot{q}_l +\dot{x}' ) -
\nabla W (q_l +x) +  \psi(q_l +x){\cal J}(\dot{q}_l +\dot{x} ) \Big|
\leq \Lambda \lambda_2^2 (r+\delta) ||x'-x||_1 e^{- 2 \lambda_2 t}.
$$
We easily get from this 
\be\label{eq:contr}
||{\cal F}(x') - {\cal F}( x) ||_1 
\leq \frac{\Lambda}{2} (r+\delta) || x'- x||_1.
\ee

By $(\ref{eq:prim}) $ and 
(\ref{eq:contr}),  to get (i) and (ii) it is enough that  $\delta$  satisfy 
$\Lambda(r+\delta)^2/4 < \delta $ as well as 
$(r+\delta) \Lambda/2 <1$.

We assume that  $r\Lambda <1/6$ (it will be  useful to 
prove the next lemma). Then it can be checked that
with $\delta=
2 r^2 \Lambda  /7 \leq r /21 $ the above inequalities are satisfied. Therefore
 equation (\ref{eq:fix}) has a solution $q_h^+ =
q_l + x$, with $||x||_1<2 \Lambda r^2 /7< r/21$.
This clearly implies the estimates of the lemma.

We must justify that $q_h^+ (\R^+) \subset B_r$. 
Using the  last estimate (but with a different $r$) and the 
uniqueness remark at the 
beginning of this proof we can get that 
$|\dot{q}_h^+ (t) + \sqrt{A} e^{-t \sqrt{A}}\beta | \leq
2 \lambda_2 \Lambda |q_h^+ (t)|^2/7 \leq 
\lambda_2 |q_h^+ (t)|/ 21$. As a consequence,
$d(|q_h^+ (t)|^2)/dt <0$ for all $t \in ( 0, + \infty) $, and we get the claim.
\end{pfn}
\\[2mm]
\begin{pfn}{\sc of lemma}\ref{lem:link1}.
We look for a solution  of (\ref{maineq}) of the form:
$$
q_d = q_h + y_l + y \  ,  \
{\rm where} \quad  
(y_l)_j  (t) = \frac{\sinh (\lambda_j t)}{\sinh (\lambda_j d)}
(\alpha - q_h (d))_j
$$ 
(we call for simplicity $ q^{+}_h = q_h $ the solution given by 
lemma \ref{lem:homo} ).
As a consequence of lemma \ref{lem:homo} and of the assumption 
$d \geq 2/\lambda_2$, we have that
\be\label{eq:yl}
|y_l (t) | \leq \frac{6}{5} r e^{-\lambda_2 (d-t)} 
\quad ; \quad 
|\dot{y}_l (t) | \leq \lambda_1
\frac{12}{5} r e^{-\lambda_2 (d-t)}. 
\ee 
We define the space:
$$
Z_2= \Big\{ g \in W^{1, \infty} [ 0 , d]  \ | \  
\sup_{t \in [0, d ]} |g(t)| e^{\lambda_2 (d - t)},
\sup_{t \in [0, d ]} |{\dot g}(t)| e^{\lambda_2 (d - t)} < +\infty 
\Big\}
$$
with norm 
$$ 
||g||_2  = \max \Big\{ \sup_{t \in [0, d ]}  |g(t)| e^{\lambda_2 (d - t)},
\sup_{t \in [0, d ]} \frac{1}{2 \lambda_2}
  |{\dot g}(t)| e^{\lambda_2 (d - t)} \Big\}.  
$$
We have to  find  a 
solution in $Z_2$ of  the fixed point problem
\be\label{eq:fix1} 
y= \wtilde{{\cal F}}(y) = \wtilde{\cal L}
\Big[ \nabla W(q_h + y_l +y) - \nabla W(q_h) -
\psi ( q_h + y_l +y ) {\cal J} ( {\dot q}_h + {\dot y}_l
+\dot{y} ) + 
\psi ( q_h  ) {\cal J} ( {\dot q}_h ) \Big], 
\ee
where $\wtilde{\cal L}$ is the linear operator which 
assigns to $h$ the unique solution $ u = \wtilde{\cal L} h $ 
of the problem
\be\label{eq:s1}
 - {\ddot u } + A u = h
\quad {\rm with } \quad u(0)= 0 \ {\rm and } \ u( d ) = 0.  
\ee
The solution $u$ of 
(\ref{eq:s1}) is  given by:
\be\label{eq:sol}
u_j (t)=  \frac{1}{ \lambda_j \sinh{(\lambda_j d) } } \left[ 
\int_{t}^d  h_j (s) \sinh( \lambda_j (d-s)) \sinh (\lambda_j t)
 ds  +  
\int_0^t  h_j (s) \sinh( \lambda_j s)\sinh ( \lambda_j (d-t)) ds \right] 
\ee
and
$$
\dot{u}_j (t)=  \frac{1}{ \sinh{(\lambda_j d)} } \left[ 
\int_{t}^d  h_j (s) \sinh( \lambda_j (d-s)) \cosh 
( \lambda_j t)ds  - 
\int_0^t  h_j (s) \sinh( \lambda_j s) \cosh 
( \lambda_j (d-t)) ds \right].
$$
It is easy to derive from these expressions the estimate
\be\label{eq:est1}
|u (t)|, \frac{1}{2\lambda_2} |\dot{u} (t)|
 \leq \frac{1}{2\lambda_2} \left[
\int_0^t |h(s) | e^{\lambda_2 (s-t)} \; ds 
+ \int_t^d |h(s) | e^{\lambda_2 (t-s)} \; ds
\right] 
\ee
As in the proof of lemma \ref{lem:homo} we have to find 
 $ \delta$ small enough
such that $\ov{\wtilde{\cal F} (B_{\delta})} \subset B_{\delta} $ 
and $\wtilde{\cal F}$ is a contraction on $B_{\delta}$. For 
$y \in B_{\delta}$ set
$$
A(t) = \nabla W (q_h + y_l + y) - 
\nabla W (q_h) \  ;  \ 
B(t) = \psi ( q_h + y_l +y ) {\cal J} ( {\dot q}_h + {\dot y}_l
+\dot{y} ) - 
\psi ( q_h  ) {\cal J} ( {\dot q}_h )
$$
We have by lemma \ref{lem:homo} and (\ref{eq:yl}):
\begin{eqnarray*}
| A(t) |& \leq & 
L_1 ( |q_h (t)| + |y_l (t) | + | y(t)| ) ( |y_l (t)|
+ |y(t)| )  \\
&\leq & L_1 \Big[ \frac{22}{21} r \Big( \frac{6}{5} r + \delta \Big) 
e^{-\lambda_2 d} +
\Big( \frac{6}{5} r + \delta \Big)^2 e^{-2\lambda_2 (d-t)} \Big],
\end{eqnarray*}
and 
\begin{eqnarray*}
| B(t)  | &\leq & L_2 \Big[
( |q_h (t)| + |y_l (t) | + | y(t)| ) ( |\dot{y}_l (t)|
+ |\dot{y}(t)| ) +
(  |y_l (t) | + | y(t)| ) |\dot{q}_h (t) | \Big] \\
&\leq & L_2 \lambda_1 \Big[ \frac{22}{7} r \Big( \frac{6}{5} r
+ \delta \Big) e^{-\lambda_2 d} + 2 \Big( \frac{6}{5} r +
\delta \Big)^2 e^{-2\lambda_2 (d-t)} \Big] .
\end{eqnarray*}
Hence
\be \label{eq:tot}
|A(t)| + |B(t)| \leq
\Lambda  \lambda_2^2 \Big[ \frac{22}{21} r \Big( \frac{6}{5} r +
\delta \Big) e^{-\lambda_2 d}  +
\Big( \frac{6}{5} r +
\delta \Big)^2   e^{-2\lambda_2 (d-t)} \Big] .
\ee
Replacing $|h(s)|$ by $|A(s)|+|B(s)|$ in
(\ref{eq:est1})  and using 
(\ref{eq:tot}), we get after easy computations :
\be \label{eq:boule}
||\wtilde{\cal F}(y)||_2 \leq \frac{\Lambda}{2} 
 \Big[ \frac{22}{21} r \Big( \frac{6}{5} r +
\delta \Big)   +
\Big( \frac{6}{5} r +
\delta \Big)^2  \Big].
\ee
We can prove in the same way that
\be \label{eq:cont}
||\wtilde{\cal F}(y)- \wtilde{\cal F} (y')||_2 \leq ||y-y'||_2 
\frac{\Lambda}{2} \Big( \frac{22}{21} r 
+ 2 \Big(\frac{6}{5} r + \delta \Big) \Big).
\ee
Using (\ref{eq:boule}) and (\ref{eq:cont}), we can see after
some elementary  calculus that, if $r\Lambda  \leq 1/10$, then
$\wtilde{\cal F}$ is a contraction on $B_{\delta}$, with
$\delta =  2  r^2 \Lambda \leq 1/5 r$.
Therefore, if $r\Lambda \leq 1/10$, we get the existence of $q_d$,
with the estimate
$$
||q_d - q_h - y_l ||_2 \leq 2  r^2 \Lambda.
$$
In particular,
$$
|q_d (t) - q_h (t)| \leq |y_l (t) | + 2r^2 \Lambda e^{-\lambda_2 
(d-t)} \ , \ 
|\dot{q}_d (t) - \dot{q}_h (t)| \leq |\dot{y}_l (t) | + 4r^2 \Lambda 
\lambda_2 e^{-\lambda_2 
(d-t)}
$$
and we get (\ref{eq:deltaQ2}) and  (\ref{eq:delta-Q2})
 by (\ref{eq:yl}). Moreover
$$
| \dot{q}_d (0) - \dot{q}_h (0) - \dot{y}_l (0)|
\leq 4  \lambda_2 r^2 \Lambda e^{-\lambda_2 d},
$$
hence 
$$
\left| ({\dot q}_d (0) - {\dot q}^{+}_h (0)) - 
 ({\dot q}_{d,L} (0) - {\dot q}^{+}_{h,L} (0)) \right|
< 4 \lambda_2  r^2 \Lambda  e^{-\lambda_2 d} + |\dot{z}(0)|,  
$$ 
where $z$ is the solution of $- {\ddot z} + Az = 0 $ 
with boundary conditions $z(0)=0$ and $z(d)= q_{h,L}(d) - q_h (d) $.
Using that  $d \geq 2/\lambda_2$,
we get by lemma \ref{lem:homo}
$$
|\dot{z}(0) | \leq 
\max_{j=1,2} \left( \frac{\lambda_j}{\sinh (\lambda_j d)} 
 \right) |q_h (d) - q_{h,L} (d) | \leq
\frac{5}{7} \lambda_2 r^2  \Lambda  e^{-\lambda_2 d}.
$$  
Estimate (\ref{eq:deltaQ}) follows.
\end{pfn}

\begin{pfn}{\sc of lemma \ } \ref{lem:nat}.
We first justufy that $w(\Theta) \in E$. Let $R_{\Theta} = Q_{\Theta} +
w (\Theta)$. By the characterisation of $w$, we have 
$S(R_{\Theta}) = \sum_{i=1}^k \alpha_i (\Theta) a_i $ for 
some $(\alpha_1 (\Theta ), \ldots , \alpha_k (\Theta ) ) \in
\R^k $. Hence, by the definition of $a_i$, outside a compact interval,
$$
-\ddot{R}_{\Theta} + \psi (R_{\Theta}) {\cal J} \dot{R}_{\Theta} 
+ A R_{\Theta} - \nabla W(R_{\Theta}) = 0. 
$$ 
Now 
$$\limsup_{|t| \to \infty} \max ( |R_{\Theta} (t)|, |\dot{R}_{\Theta} (t)|
/ \lambda_2 ) = \limsup_{|t| \to \infty} 
\max ( |w(\Theta) (t)|, |\dot{w}({\Theta}) (t)|
/ \lambda_2 ) < 
\min(2/\Lambda , \rho_0 ),$$
by lemma \ref{lem:const}. Arguing as for the proof of lemma \ref{lem:one},
one can derive that $\lim_{|t| \to \infty} |R_{\Theta}(t)| +
|\dot{R}_{\Theta} (t) | =0 $, which implies, by the properties
of the equation near the equilibrium, that $R_{\Theta} \in X$. 
Hence $w(\Theta ) \in X \subset E$.

Assume that $\ov{\Th}$ is a critical point of $f( Q_{\Th} + w_{\Th})$
(equivalently $(\ov{d}_1, \ldots , \ov{d}_k ) $ is a critical point of $g(d)$).
Then there results that:
\be\label{eq:tang}
\Big( S( Q_{\ov{\Th}} + w(\ov{\Th})), \frac{\partial}{\partial \th_j} 
( Q_{\ov{\Th}} + w(\ov{\Th}) ) \Big) = 0 \quad {\rm for } \quad 
j = 1, \ldots , k .
\ee  
Hence  we have that 
\be\label{eq:comb}
\sum_{i=1}^k \al_i (\ov{\Th})(a_i , \frac{\partial}{\partial \theta_j}
(Q_{\ov{\Teta}} + w (\ov{\Teta}) ) = \sum_{i=1}^k \al_i (\ov{\Th}) \alpha  
\int_{J_i} \dot{Q}^i (t) \cdot  \frac{\partial}{\partial \theta_j}
(Q_{\ov{\Teta}} + w (\ov{\Teta}) )(t) \ dt = 0.
\ee
Letting 
$
b_{i,j} (\Th ) = \int_{J_i}  \dot{Q}^i (t) \cdot  
\frac{\partial}{\partial \theta_j}
(Q_{\Teta} + w (\Teta) )(t) \ dt,
$
(\ref{eq:comb}) yields 
$\sum_{i=1}^k \al_i (\ov{\Th}) b_{i,j} (\ov{\Th}) = 0.$
Now we show that 
\be\label{eq:bij}
b_{i,j}(\Th) = \delta_{i,j}
 \displaystyle \int_{J_j} -|\dot{Q}^j (t) |^2 - \dot{Q}^j (t)
\cdot \dot{w} (\Teta)(t) dt 
\ee 
It appears that $ (\partial  Q_{\Teta |J_i} / \partial \teta_j) = $ 
$- \delta_{i,j} \dot{Q^i}_{|J_i} \in $ 
$W^{1, \infty} (J_i)$. 
Hence, since
 $(Q_{\Teta} + w(\Teta))$ is a $C^1$ function of $\Theta$, 
$\partial w(\Teta)_{|J_i} / \partial \teta_j$ too is well defined in 
$W^{1, \infty} (J_i)$. 
We know that 
\be \label{eq:0qi}
0=\int_{J_i} \dot{Q}^i (t) \cdot w(\Teta)(t) \; dt =
\int_{J_i - \Teta_i} \dot{Q}^i (t+ \teta_i) \cdot w(\Teta)
(t+ \teta_i) \; dt. 
\ee 
Now $J_i - \teta_i $ and $\dot{Q}^i (\cdot + \teta_i) $
do not depend on $\Teta$. Hence, deriving (\ref{eq:0qi}) with
respect to $\teta_j$ we get 
$$ 
\int_{J_i} \dot{Q}^i (t) \cdot \frac{\partial}{\partial \teta_j} 
w(\Th) (t) \ dt 
 = - \delta_{i,j}
 \int_{J_i} \dot{Q}^i (t)  \cdot \dot{w} ({\Teta}) (t) dt .
$$   
So  (\ref{eq:bij}) is proved.
Therefore,  if $ \ov{\Teta} $ is
a critical point of $f(Q_{\Teta} + w(\Teta))$, then 
$\alpha_j ( \ov{\Th} )  b_{j,j} (\ov{\Th})  =0$ 
for all $j$, and it is enough to 
check that $b_{j,j} (\ov{\Th}) \neq 0 $ to conclude.
But, if for example $ Q^i (t) = \ov{q}( t - \ov{\th}_i )$ then  
$$|b_{j,j}(\ov{\Th}) | \geq \int_{{[ \ov{t} - \ov{\tau}, \ov{t} + \ov{\tau}]}}
\frac{3}{4} |\dot{\ov{q}} (\ov{t})| \Big( \frac{3}{4}
|\dot{\ov{q}} (\ov{t} )| - \lambda_2 ||w(\ov{\Th})||  \Big) dt >0
$$   
because, for $\ov{d} \geq \ov{D}$, 
$||w(\ov{\Teta})|| \leq r/2$ and, by lemma {lem:homo} 
$ |\dot{\ov{q}}(\ov{t} )|
\geq 21 \lambda_2 r / 22$. 
\end{pfn}

\begin{pfn}{\sc of lemma \ }\ref{lem:en}.
Consider
\be\label{eq:act}
e(\bt, \al, d)=
\frac{1}{2} \int_0^d {\dot q}^2_d (t) + A q_d (t) \cdot q_d (t) dt 
- \int_0^d {\dot q}_d \cdot v(q_d) dt - \int_0^d W(q_d(t)) dt.
\ee
We perform the change of variables $t=sd$ and set 
$$ 
\forall s \in [0,1] \  Q_d (s) = q_d (sd ) .
$$
The function $Q_d (s)$ satisfies the equation:
\be\label{eq:tras}
- {\ddot Q_d} + d^2 A Q_d + d \psi( Q_d) {\cal J} {\dot Q}_d
= d^2 \nabla W ( Q_d )  \quad 
{\rm with } \quad Q_d(0)=\bt \ {\rm and } \ Q_d(1) = \al.
\ee
After the change of variables in (\ref{eq:act}) we will have
$$
e(\bt, \al , d) =
\frac{1}{2} \int_0^1 \frac{1}{d} {\dot Q}^2_d (s) + 
d A Q_d (s) \cdot Q_d (s)ds - \int_0^1 d W (Q_d(s)) ds - 
\int_{0}^{1} {\dot Q}_d(s) \cdot v (Q_d (s))ds   
$$
We now take the derivative of $e(\bt, \al, d )$ with respect to
$d$. It is given by
\be\label{eq:es}
\frac{ \partial e}{\partial d} (\bt, \al , d) =
\int_0^1 - \frac{1}{2 d^2} {\dot Q}^2_d (s) +
\frac{1}{d} {\dot Q}_d (s) \cdot \partial_d {\dot Q}_d (s)  +
\frac{1}{2} A Q_d (s) \cdot Q_d (s) +
d A Q_d (s) \cdot {\partial}_d Q_d (s) +
$$
$$
- \int_0^1  W (Q_d(s)) -
d \nabla W (Q_d (s)) \partial_d Q_d (s) ds - 
\frac{\partial}{\partial d } \int_0^1  {\dot Q}_d \cdot v( Q_d ) ds 
\ee
After an integration by parts in \ref{eq:es} (of the term 
$ {\dot Q}_d (s) \cdot \partial_d {\dot Q}_d (s) / d $)
using  that $ \partial_d Q_d (0)=0 \quad  {\rm and } \quad 
\partial_d Q_d(1) = 0 $
and putting (\ref{eq:tras}) into (\ref{eq:es}) we have that
\be\label{eq:es1}
\frac{ \partial e}{\partial d} (\bt, \al , d) =
\int_0^1  - \frac{1}{2 d^2} {\dot Q}^2_d  +
\frac{1}{2} A Q_d  \cdot Q_d  -  W (Q_d)  ds 
- \int_0^1 \psi(Q_d) \partial_d Q_d {\cal J} {\dot Q}_d - 
\frac{\partial}{\partial d } \int_0^1 v( Q_d) {\dot Q}_d ds 
\ee
Using  that 
$$
\int_0^1 \psi(Q_d) \partial_d Q_d  {\cal J} {\dot Q}_d + 
\frac{\partial}{\partial d } \int_0^1 v( Q_d) {\dot Q}_d ds =0. 
$$
and making the change of variable $t=sd$ in (\ref{eq:es1}), we get 
$$
\frac{ \partial e}{\partial d} (\bt, \al , d) =
- \int_0^d \left[ \frac{1}{2} ( {\dot q}^2_d (t) -
A q_d (t) \cdot q_d (t) ) + W (q_d (t)) \right] \frac{1}{d}  dt.   
$$
Since the integrant is nothing but the energy
$ ( {\dot q}^2_d (t) -
A q_d (t) \cdot q_d (t) ) / 2  + W (q_d (t)) =  {\cal E}(d) $
which is constant we finally get that $ \displaystyle
\frac{ \partial e}{\partial d} (\bt, \al , d) = - {\cal E}(d)d \frac{1}{d} =
- {\cal E}(d).$ 
\end{pfn}

{\it Massimiliano Berti, Scuola Normale Superiore, Pza dei Cavalieri
7, 56100, Pisa, berti@cibs.sns.it}. 
\\[2mm]
{\it Philippe Bolle, Dep. of Math. Sc., University of Bath,
Bath BA2 7AY, maspb@maths.bath.ac.uk}.
\end{document}